\documentclass[11pt]{article}
\usepackage[normalem]{ulem}
\usepackage[left=1.2in,top=1.5in,right=1.2in,bottom=1.5in]{geometry}
\usepackage{amsmath}
\usepackage{latexsym, amsbsy, amssymb, amsmath, amsthm, amsfonts}
\usepackage{graphicx}
\usepackage{natbib}
\usepackage{fancybox,fancyheadings}
\usepackage{psfrag}
\usepackage{float}
\usepackage{titlesec}
\usepackage[notlof,notlot,nottoc]{tocbibind}
\usepackage{cases}
\usepackage{sectsty}
\usepackage{chngcntr}
\usepackage{undertilde}
\usepackage{enumerate}
\usepackage{enumitem}
\usepackage{tabulary}
\usepackage{caption}
\usepackage[toc,titletoc]{appendix}
\usepackage{multirow}
\usepackage{multicol}
\usepackage{rotating}
\usepackage{longtable}
\usepackage{fancyhdr}
\RequirePackage[colorlinks,breaklinks]{hyperref} \hypersetup{
    colorlinks = true,
    linkcolor = blue,
    anchorcolor = red,
    citecolor = blue,
    filecolor = red,
    urlcolor = blue
}

\newcommand{\bm}[1]{\mbox{\boldmath{$#1$}}}
\newcommand{\matern}[0]{$\rm{Mat\acute{e}rn}\ $}
\newcommand{\Bmatern}[0]{$\bf{Mat\acute{e}rn}\ $}
\newcommand{\ScriptSmallRM}[1]{\mbox{\tiny$\rm{#1}$}}
\newcommand{\ScriptSmallCap}[1]{\textsc{#1}}
\newcommand{\ScriptSmall}[1]{\mbox{\tiny${#1}$}}
\newcommand{\LargerCdot}[1]{\raisebox{-0.25ex}{\scalebox{#1}{$\cdot$}}}

\newenvironment{mydescription}[1]{\setdescription{leftmargin = #1, labelindent = 0pt}\begin{description}}{\end{description}}
\allowdisplaybreaks
\numberwithin{equation}{section}

\title{Functional Principal Components Analysis of Spatially Correlated Data}

\author{Chong Liu, Surajit Ray and Giles Hooker}

\begin{document}

\maketitle

\begin{center}
ABSTRACT
\end{center}
\small
This paper focuses on the analysis of spatially correlated functional data. The between-curve correlation is modeled by correlating functional principal component scores of the functional data. We propose a Spatial Principal Analysis by Conditional Expectation framework to explicitly estimate spatial correlations and reconstruct individual curves. This approach works even when the observed data per curve are sparse. Assuming spatial stationarity, empirical spatial correlations are calculated as the ratio of eigenvalues of the smoothed covariance surface $Cov(X_i(s),X_i(t))$ and cross-covariance surface $Cov(X_i(s), X_j(t))$ at locations indexed by $i$ and $j$. Then a anisotropy \matern spatial correlation model is fit to empirical correlations. Finally, principal component scores are estimated to reconstruct the sparsely observed curves. This framework can naturally accommodate arbitrary covariance structures, but there is an enormous reduction in computation if one can assume the separability of temporal and spatial components. We propose hypothesis tests to examine the separability as well as the isotropy effect of spatial correlation. Simulation studies and applications of empirical data show improvements in the curve reconstruction using our framework over the method where curves are assumed to be independent. In addition, we show that the asymptotic properties of estimates in uncorrelated case still hold in our case if 'mild' spatial correlation is assumed.

\normalsize
\fancypagestyle{plain}{
\fancyhf{}
\renewcommand{\headrulewidth}{0pt}
\renewcommand{\footrulewidth}{0pt}}

\pagestyle{fancy}
\lhead{}
\chead{\thepage}
\rhead{}
\lfoot{}
\cfoot{}
\rfoot{}
\renewcommand{\headrulewidth}{0pt}
\renewcommand{\footrulewidth}{0pt}
\pagenumbering{arabic}
\setcounter{page}{1}
\section{\bf{\large Introduction}}\label{section: introduction}
	Functional data analysis (FDA) focuses on data that are infinite-dimensional, such as curves, shapes and images. Generically, functional data are measured over a continum across multiple subjects. In practice, many data such as growth curves of different people, gene expression profiles, vegetation index across multiple locations, vertical profiles of atmospheric radiation recorded at different times, etc. could naturally be modeled by FDA framework.
		
	Functional data are usually modeled as noise corrupted observations from a collection of trajectories that are assumed to be realizations of a smooth random function of time $X(t)$, with unknown mean shape $\mu(t)$ and covariance function $Cov(X(s),X(t)) = G(s,t)$. The functional principal components (fPCs) which are the eigenfunctions of the kernel $G(s,t)$ provide a comprehensive basis for representing the data and hence are very useful in problems related to model building and prediction of functional data.
	
	Let ${\phi_k(t),\ k = 1,2,\cdots,K}$ and ${\lambda_k,\ k = 1,2,\cdots,K}$ be the first $K$ eigenfunctions and eigenvalues of $G(s,t)$. Then 
	$$
		X_i(t)\approx\sum_{k=1}^K\xi_{ik}\phi_k(t)
	$$
	where $\xi_{ik}$ are fPC scores which have mean zero and variance $\lambda_k$. According to this model, all curves share the same mode of variations, $\phi_k(t)$, around the common mean process $\mu(t)$.
	
	A majority of previous work in FDA assume that the realizations of the underlying smooth random function are independent. There exists an extensive literature on functional principal components analysis (fPCA) for this case.  For data observed at irregular grids, \citet{YaoMullerClifford2003} and \citet{YaoLee2005} used local linear smoother to estimate the covariance kernel and integration method to compute fPC scores. However, the integration approximates poorly with sparse data. \citet{JamesSugar2003} proposed B-splines to model the individual curves through mixed effects model where fPC scores are treated as random effects. For sparsely observed data, \citet{YaoMullerWang2005} proposed a framework called ``PACE'' which stands for Principal Analysis of Condition Expectation. In  PACE, fPC scores were estimated by their expectation conditioning on available observations across all trajectories. To estimate fPCs: a  system of orthogonal functions, \citet{PengPaul2009} proposed a restricted maximum likelihood method based on a Newton-Raphson procedure on the Stiefel manifold. \citet{HallMullerWang2006} and \citet{LiHsing2010} gave weak and strong uniform convergence rate of the local linear smoother of the mean and covariance, and the rate of derived fPC estimates.

	 The PACE approach works by efficiently extracting the information on $\phi_k(t)$ and $\mu(t)$ even when only a few observations are made on each curve as long as the pooled time points from all curves are sufficiently dense. Nevertheless, PACE is limited by its assumption of independent curves. In reality, observations from different subjects are correlated. For example, it is expected that expression profiles of genes involved in the same biological processes are correlated; and vegetation indices of the same land cover class at neighboring locations are likely to be correlated. 
		
	There has been some recent work on correlated functional data. \citet{li2007nonparametric} proposed a kernel based nonparametric method to estimate correlation among functions where observations are sampled at regular temporal grids and smoothing is performed across different spatial distances. Moreover, it was assumed in their work that the covariance between two observations can be factored as the product of temporal covariance and spatial correlation, which is referred to as separable covariance. \citet{Paul2010Manu} discussed a nonparametric method similar to PACE to estimate fPCs and proved that the $L^2$ risk of their estimator achieves optimal nonparametric rate under mild correlation regime when the number of observations per curve is bounded. \citet{Zhou2010Mixed} presented a mixed effect model to estimate correlation structure, which accommodates both separable and non-separable structures. 
		
	In this paper, we develop a new framework which we call SPACE (Spatial PACE for modeling correlated functional data. In SPACE, we explicitly model the spatial correlation among curves and extend local linear smoothing techniques in PACE to the case of correlated functional data. Our method differs from \citet{li2007nonparametric} in that sparsely and irregularly observed data can be modeled and it is not necessary to assume separable correlation structure. In fact, based on our SPACE framework, we proposed hypothesis tests to examine whether or not correlation structure presented by data is separable or not.    
	
	Specifically, we model the correlation of fPC scores $s_{ik}$ across curves by anisotropiv \matern family. In the anisotropy \matern correlation model \citep{haskard2007anisotropic}, we rotate and stretch the axis such that equal correlation contour is a tilted ellipse to accommodate anisotropy effect which often arises in geoscience data. In our model, anisotropy \matern correlation is governed by 4 parameters: $\alpha,\delta,\kappa,\phi$ where $\alpha$ controls the axis rotation angle and $\delta$ specifies the amount of axis stretch. SPACE identifies a list of neighborhood structures and applies local linear smoother to estimate a cross-covariance surface for each spatial separation vector. An example of neighborhood structure could be all pairs of locations which are separated by distance of one unit and are positioned from southwest to northeast. In particular, SPACE estimates a cross-covariance surface by smoothing empirical covariances observed at those locations. Next, empirical spatial correlations are estimated based on the eigenvalues of those cross-covariance surfaces. Then, anisotropy \matern parameters are estimated from the empirical spatial correlations. SPACE directly plugs in the fitted spatial correlation model into curve reconstruction to improve the reconstruction performance relative to PACE where no spatial correlation is modeled.
	
	We demonstrate SPACE methodology using simulated functional data and Harvard Forest vegetation index discussed in \citet{liu2012functional}. In simulation studies, we first examine the estimation of SPACE model components. Then we perform the hypothesis tests of separability and isotropy effect. We show that curve reconstruction performance is improved using SPACE over PACE. Also, hypothesis tests demonstrate reasonable sizes and powers. Moreover, we construct semi-empirical data by randomly removing observations to achieve sparseness in vegetation index at Harvard Forest. Then it is shown that SPACE restores vegetation index trajectories with less errors than PACE. 

	The rest of the paper is organized as follows. Section \ref{section: correlated functional data model} describes the spatially correlated functional data model. Section \ref{section: space methodology} describes the SPACE framework and model selections associated with it. Then we summarize the consistency results of SPACE estimates in Section \ref{section: consistency} and defer more detailed discussions to Appendix \ref{appendix: consistency}. Next, we propose hypothesis tests based on SPACE model in Section \ref{section: separability and isotropy tests}. Section \ref{section: simulation} describes simulation studies on model estimations, followed by Section \ref{section: Harvard Forest} which presents curve construction analysis on Harvard Forest data. In the end, conclusion and comments are given in Section \ref{section: conclusion}. 
	
	\section{\bf{\large Correlated Functional Data Model}}\label{section: correlated functional data model}
		In this section, we describe how we incorporate spatial correlation into functional data and introduce the \matern class which we use to model spatial correlation.
		\subsection{\bf{\large Data Generating Process}}\label{subsection: data generating process}		
		We start by assuming that data are collected across $N$ spatial locations. For location $i$, a number of $n_i$ noise-corrupted points are sampled from a random trajectory $X_i(t)$, denoted by $Y_i(t_j),\ j = 1,2,\cdots,n_i$. These observations can be expressed by an additive error model as the following,
			\begin{equation}\label{added error model}
				Y_i(t) = X_i(t) + \epsilon_i(t).
			\end{equation}
			Measurement errors $\{\epsilon_i(t_j)\}_{i=1\ j=1}^{N\ \ \   n_i}$ are assumed to be iid with variance $\sigma^2$ across locations and sampling times. The random function $X_i(t)$ is the $i$th realization of an underlying random function $X(t)$ which is assumed to be smooth and square integrable on a bounded and closed time interval $\cal T$. Note that we refer to the argument of function as time without loss of generality. The mean and covariance functions of $X(t)$ are unknown and denoted by $\mu(t) = E(X(t))$ and $G(s,t) = Cov(X(s), X(t))$. By the Karhunen-Lo$\rm{\grave{e}}$ve theorem, under suitable regularity conditions, there exists an eigen-decomposition of the covariance kernel $G(s,t)$ such that 
			\begin{equation}
				G(s,t) = \sum_{k = 1}^{\infty}\lambda_k\phi_k(s)\phi_k(t),\ t,s\in{\cal T}
			\end{equation} 
			where $\{\phi_k(t)\}_{k=1}^{\infty}$ are orthogonal functions in the $L^2$ sense which we also call functional principal components (fPC), and $\{\lambda_k\}_{k=1}^{\infty}$ are associated non-increasing eigenvalues. Then, each realization $X_i(t)$ has the following expansion,
			\begin{equation}
				X_i(t) = \mu(t) + \sum_{k = 1}^{\infty}\xi_{ik}\phi_i(t),\ i = 1,2,\cdots,N
			\end{equation} 
			where for given $i$, $\xi_{ik}$'s are uncorrelated fPC scores with variance $\lambda_k$. Usually, a finite number of eigenfunctions are chosen to achieve reasonable approximation. Then,
			\begin{equation}\label{approx expansion}
				X_i(t) \approx \mu(t) + \sum_{k = 1}^{K}\xi_{ik}\phi_i(t),\ i = 1,2,\cdots,N	
			\end{equation}
			In classical functional data model, $X_i(t)$'s are independent across $i$ and thus $cor(\xi_{ik}, \xi_{jk}) = 0$ for any pair of different curves $\{i,j\}$ and for any given fPC index $k$. However, in many applications, explicit modeling and estimation of the spatial correlation is desired and can provide insights into subsequent analysis. To build in correlation among curves, we assume $\xi_{ik}$'s are correlated across $i$ for each $k$. One could specify full correlation structure among $\xi_{ik}$'s by allowing non-zero covariance between scores of different fPCs, e.g. $Cov(\xi_{ip}, \xi_{jq})\neq 0$. Though the full structure is very flexible, it is subject to the risk of overfitting and thus its estimation can be intractable. To achieve parsimony, we assume the following
			\begin{equation}\label{raw spatial correlation}
  			Cov(\xi_{ip}, \xi_{jq}) = \left\{\begin{array}{ll}
    																\rho_{ij}(k)\lambda_k,\ \ & \text{if $p = q = k$},\\\\
    																0, & \text{otherwise},
  						 										\end{array}\right.
			\end{equation} 
			where $\rho_{ij}(k)$ measures the correlation between $k$th fPC scores at curve $i$ and $j$. Denoting $\bm\xi_i = (\xi_{i1},\xi_{i2},\cdots,\xi_{iK})^T$, $\bm\phi(t) = (\phi_1(t),\phi_2(t),\cdots,\phi_K(t))^T$ and retaining the first $K$ eigenfunctions as in \eqref{approx expansion}, then the covariance between $X_i(s)$ and $X_j(t)$ can be expressed as 
			\begin{eqnarray}\label{cross covariance}
				Cov(X_i(s), X_j(t))	&=& \bm\phi(s)^Tcov(\bm\xi_i\bm\xi_j^T)\bm\phi(t)\\ 
				&=& \bm\phi(s)^T{\textrm{diag}}\left(\rho_{ij}(1)\lambda_1,\rho_{ij}(2)\lambda_2,\cdots,\rho_{ij}(K)\lambda_K\right)\bm\phi(t),
			\end{eqnarray}
			Note that $Cov(\bm\xi_i\bm\xi_j^T)$ is a diagonal matrix. Hence, columns of $\bm\phi(t)$ are comprised of eigenvectors of $Cov(X_i(s), X_j(t))$. Note the diagonal elements are not necessarily sorted by value. If we further assume the between-curve correlation $\rho_{ij}(k)$ does not depend on $k$, i.e. $\rho_{ij}(k) = \rho_{ij}$, then $Cov(\bm\xi_i\bm\xi_j^T) = \rho_{ij}\rm{diag}(\lambda_1,\cdots,\lambda_K)$. In this case the covariance can be further simplified as 
			\begin{equation}
				Cov(X_i(s), X_j(t)) = \rho_{ij}\bm\phi(s)^T{\rm diag}(\lambda_1,\cdots,\lambda_K)\bm\phi(t) = \rho_{ij}Cov(X(s), X(t)).\label{cor.separable}
			\end{equation} 
			If the covariance between $X_i(s)$ and $X_j(t)$ can be decomposed into a product of spatial and temporal components as in \eqref{cor.separable}, we refer to this covariance structure as separable. Separable covariance structure of the noiseless processes assumes that the correlation across curves and across time are independent of each other. One example of this type of processes is the weed growth data studied in \citet{banerjee2006} where curves are weed growth profiles at different locations in the agricultural field. Separable covariance is a convenient assumption which makes estimation easier. It is a restrictive assumption and in order to examine whether this assumption is justifiable, we propose a hypothesis test which is described in Section \ref{section: separability and isotropy tests}. Irrespective of the separability of covariance, spatial correlations among curves are reflected through the correlation structure of fPC scores $\bm\xi$. For each fPC index $k$, the associated fPC scores at different locations can be viewed as a spatial random process.
			
			\subsection{\bf{\large\Bmatern Class}}\label{subsection: matern class}	
			We choose the \matern class for modeling spatial correlation. It is a widely used class as a parametric model in geoscience. The \matern family is attractive due to its flexibility. Specifically, the \matern correlation between two observations at locations separated by distance $d > 0$ is given by
			\begin{equation}
				\rho(d;\zeta,\nu) = \dfrac{1}{2^{\nu-1}\Gamma(\nu)}\left(\dfrac{d}{\zeta}\right)^{\nu}K_{\nu}\left(\dfrac{d}{\zeta}\right),
			\end{equation}  
			where $K_{\nu}(\cdot)$ is the modified Bessel function of the third kind of order $\nu > 0$ described in \citet{abramowitz1970handbook}. This class is governed by two parameters, a range parameter $\zeta > 0$ which rescales the distance argument and a smoothness parameter $\nu > 0$ that controls the smoothness of the underlying process. 
			
\begin{figure}[!h]
				\begin{minipage}[t]{0.48\linewidth}
    			\centering
    			\includegraphics[width = 0.95\linewidth]
    			{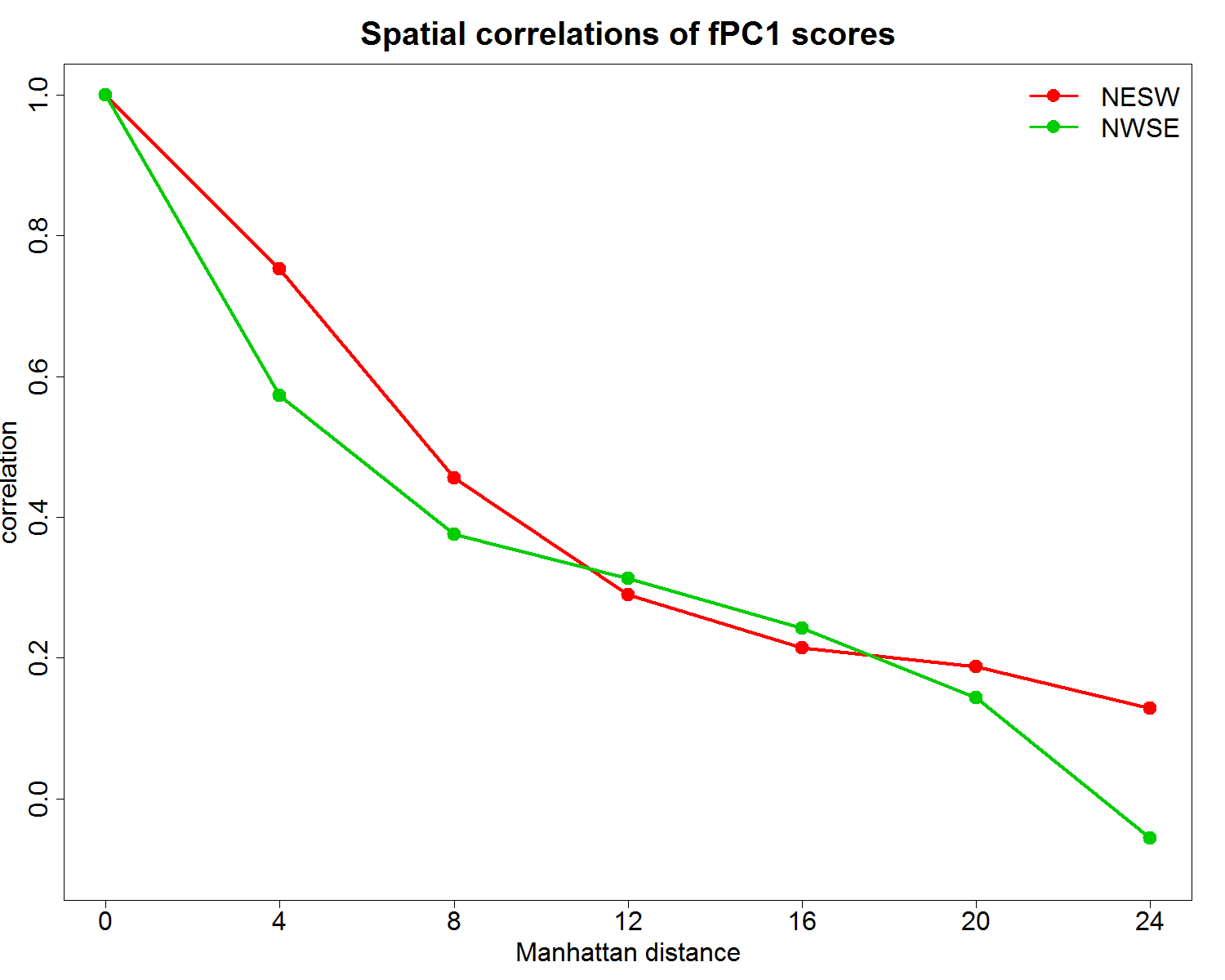}
  			\end{minipage}
  			\hspace{5pt}
  			\begin{minipage}[t]{0.48\linewidth}
    			\centering
    			\includegraphics[width = 0.95\linewidth]
      		{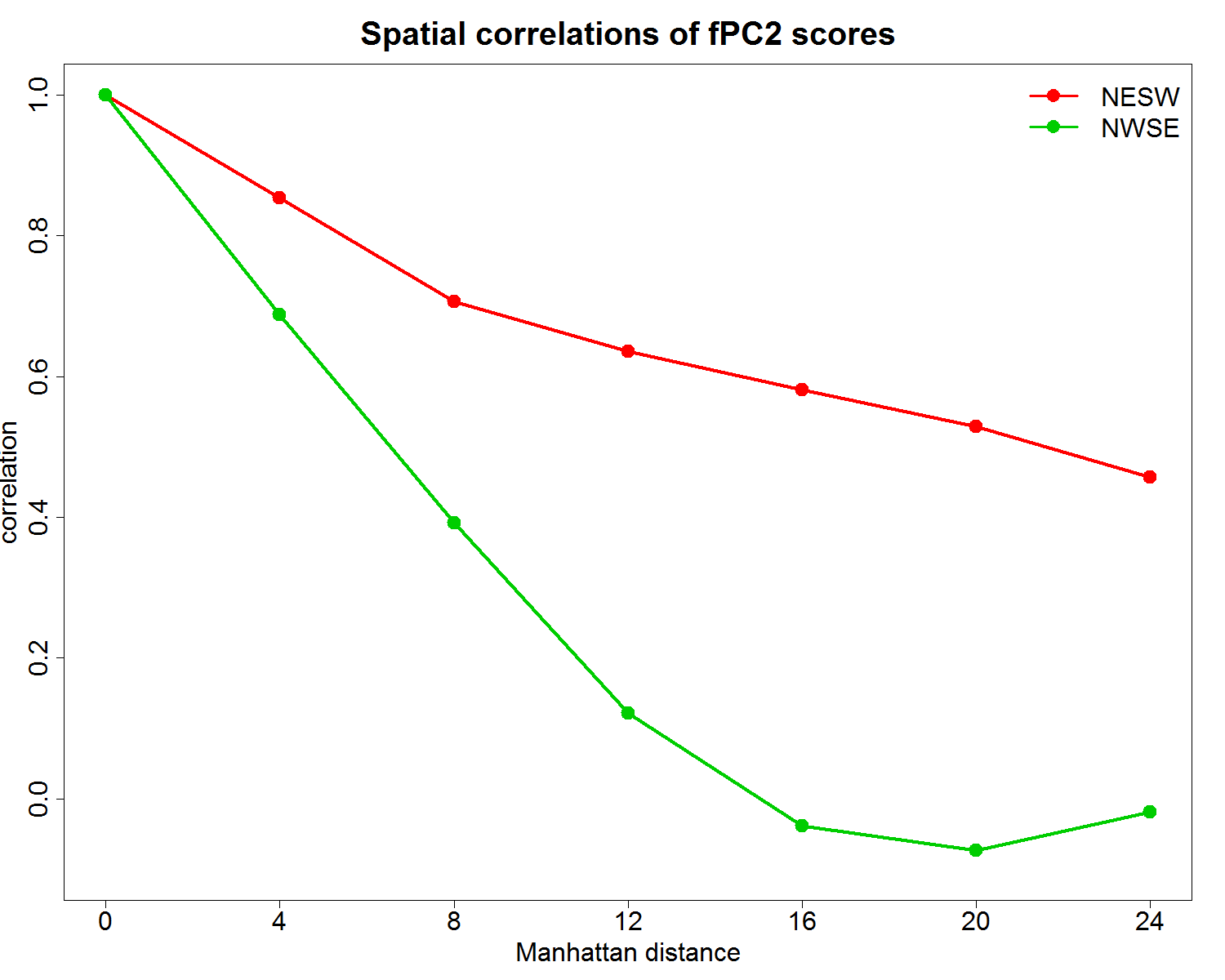}
  			\end{minipage}
    		\caption[Spatial correlations of leading fPC scores computed in directions of northeast-southwest and northwest-southeast separately at Harvard Forest]{\footnotesize{Spatial correlations of leading fPC scores are computed in directions of northeast-southwest and northwest-southeast separately at Harvard Forest. Consider a two-dimensional coordinates of integers. We calculate the correlation of fPC scores at locations which are positioned along either northwest-southeast or northeast-southwest diagonals respectively, and are separated by increasing Manhattan distances. Green lines indicate the correlations along the northwest- southeast diagonal and red lines indicate the correlations along the  northeast-southwest diagonal. The first two leading fPCs are illustrated.}}
    			\label{fig: HarvardfPCscore}
  		\end{figure}

			The \matern class is by itself isotropic which has contours of constant correlation that are circular in two-dimensional applications. However, isotropy is a strong assumption and thus limit the model flexibility in some applications. \citet{liu2012functional} showed that fPC scores present directional patterns which might be associated with geographical features. Specifically, the authors look at vegetation index series at Harvard Forest in Massachusetts across a 25 $\times$ 25 grid over 6 years. Using the same data, we calculate the correlation between fPC scores at locations separated by 45 degree to the northeast and to the northwest respectively. Figure \ref{fig: HarvardfPCscore} suggests the anisotropy effect in the second fPC scores.	
			Geometric anisotropy can be easily incorporated into the \matern class by applying a linear transformation to the spatial coordinates. To this end, two additional parameters are required: an anisotropy angle $\alpha$ which determines how much the axes rotate clockwise, and an anisotropy ratio $\delta$ specifying how much one axis is stretched or shrunk relative to the other. Let $(x_1, y_1)$ and $(x_2, y_2)$ be coordinates of two locations and denote the spatial separation vector between these two locations by $\bm\Delta = (\Delta x, \Delta y)^T = (x_2 - x_1, y_2 - y_1)^T$. The isotropy \matern correlation is computed typically based on Euclidean distance and we have $\rho(d;\zeta,\nu) = \rho(\sqrt{\bm\Delta^T\bm\Delta};\zeta,\nu)$. To introduce anisotropy correlation, a non-Euclidean distance could be applied to variable $d$ through linear transformations on coordinates. Specifically, we form the new spatial separation vector
			\begin{equation}
				\bm\Delta^* = \left(\begin{array}{c}{\Delta x}^*\\{\Delta y}^*\end{array}\right) = \left(\begin{array}{cc}\sqrt{\delta} & 0\\0 & \dfrac{1}{\sqrt{\delta}}\end{array}\right)\left(\begin{array}{cc}\cos\alpha & 	
				\sin\alpha\\-\sin\alpha & \cos\alpha\end{array}\right)\left(\begin{array}{c}{\Delta x}\\{\Delta y}\end{array}\right) = {\bf S}{\bf R}\bm\Delta	
			\end{equation}
where we write the rotation and rescaling matrix by ${\bf R}$ and ${\bf S}$ respectively. Define the non-Euclidean distance function as $d^*(\bm\Delta,\alpha,\delta) = \sqrt{\bm\Delta^{*T}\bm\Delta^*} = \sqrt{\bm\Delta^T{\bf R}^T{\bf S}^2{\bf R}\bm\Delta}$. Hence, the anisotropy correlation is computed as 
			\begin{equation}\label{anisotropy matern correlation}
				\rho^*(\bm\Delta;\alpha,\delta,\zeta,\nu) = \rho(d^*(\bm\Delta,\alpha,\delta);\zeta,\nu) = \rho(\sqrt{\bm\Delta^T{\bf R}^T{\bf S}^2{\bf R}\bm\Delta};\zeta,\nu)
			\end{equation} 		
Note $\rho^*(\bm\Delta) = \rho^*(-\bm\Delta)$. Let $\bm\Delta_{ij}$ be the spatial separation vector of locations between curve $i$ and $j$. Then we model the covariance structure of fPC scores described in \eqref{raw spatial correlation} as
			\begin{equation}\label{spatial correlation}
  			Cov(\xi_{ip}, \xi_{jq}) = \left\{\begin{array}{ll}
    																\rho^*(\bm\Delta_{ij};\alpha_k,\delta_k,\zeta_k,\nu_k)\lambda_k = \rho^*_k(\bm\Delta_{ij})\lambda_k,\ \ & \text{if $p = q = k$},\\\\
    																0, & \text{otherwise}.
  						 										\end{array}\right.
			\end{equation} 	
			Note that the above parametrization is not identifiable. Firstly, $\alpha$ and $\alpha$ + $\pi$ always give the same correlation. Secondly, the pair of any $\alpha$ and $\delta$ gives the exact same correlation with the combination of $\alpha + \pi/2$ and $1/{\delta}$. We remove the non-uniqueness of model by adding additional constraints on the ranges of parameters.
			
	\section{\bf{\large SPACE Methodology}}\label{section: space methodology}
			SPACE methodology extends the PACE methodology introduced by \citet{YaoMullerWang2005} and the methodology discussed in \citet{li2007nonparametric}. Among the components of SPACE, the mean function $\mu(t)$ and the measurement error variance $\sigma^2$ in \eqref{added error model} are estimated by the same method used in \citet{YaoMullerWang2005}. In this section, we introduce the estimation of the other key components of SPACE: the cross-covariance surface $Cov(X_i(s), X_j(t))$ for any location pair $(i,j)$ in \eqref{cross covariance}  and the anisotropy \matern model among fPC scores in \eqref{spatial correlation}. We will also describe methods of curve reconstruction and model selection.
			
		\subsection{Cross-Covariance Surface}
			Let us define $G_{ij}(s,t) = Cov(X_i(s), X_j(t))$ as the cross covariance surface between location $i$ and $j$. Let $\hat\mu(t)$ be an estimated mean function, and then let $D_{ij}(t_{ik},t_{jl}) = (Y_i(t_{ik}) - \hat\mu(t_{ik}))(Y_j(t_{jl}) - \hat\mu(t_{jl}))$ be the raw cross covariance based on observations from curve $i$ at time $t_{ik}$ and curve $j$ at time $t_{jl}$. We estimate $G_{ij}(s,t)$ by smoothing $D_{ij}(t_{ik},t_{jl})$ using a local linear smoother. 
			
			In this work, we assume the second order spatial stationarity of the fPC score process.	Define $N(\bm\Delta) = \{(i,j),\ s.t.\ \bm\Delta_{ij} = (\Delta_x, \Delta_y)\ {\textrm{or}}\ \bm\Delta_{ij} = (-\Delta_x, -\Delta_y)\}$ to be the collection of location pairs that share the same spatial correlation structure. Then, all location pairs that belong to $N(\bm\Delta)$ are associated with the same unique covariance surface which we write as $G_{\ScriptSmallCap{\bm\Delta}}(s,t)$. As a result, all raw covariances constructed based on locations in $N(\bm\Delta)$ can be pooled together to estimate $G_{\ScriptSmallCap{\bm\Delta}}(s,t)$. In addition, we note 
			\begin{equation}\label{expected raw covariance}
				E(D_{ij}(t_{ik},t_{jl}))\approx G_{ij}(t_{ik},t_{jl}) + \delta(i = j, s = t)\sigma^2,
			\end{equation}
			where $\delta(i = j, s = t)$ is 1 if $i = j$ and $s = t$, and 0 otherwise. If $i = j$, the problem reduces to the estimation of covariance surface and we apply the same treatment described in \citet{YaoMullerWang2005} to deal with the extra $\sigma^2$ on the diagonal. For $i\neq j$ and a given spatial separation vector $\bm\Delta$, the local linear smoother of the cross covariance surface $G_{\ScriptSmallCap{\bm\Delta}}(s,t)$ is derived by minimizing
			\begin{equation}\label{2d cross cov smoother}
				\sum_{(i,j) \in N(\ScriptSmallCap{\bm\Delta})}\sum_{k = 1}^{n_i}\sum_{l = 1}^{n_j}\kappa_2(\dfrac{t_{ik} - s}{h_G},\dfrac{t_{jl} - t}{h_G})(D_{ij}(t_{ik},t_{jl}) - \beta_0 - \beta_1(s - t_{ik}) - \beta_2(t - t_{jl}))^2,
			\end{equation}
			with respect to $\beta_0,\beta_1$ and $\beta_2$. $\kappa_2$ is the two-dimensional Gaussian kernel. Let $\hat\beta_0$, $\hat\beta_1$ and $\hat\beta_2$ be minimizers of \eqref{2d cross cov smoother}. Then $\widehat G_{\ScriptSmallCap{\bm\Delta}}(s,t) = \hat\beta_0$. For computation, we evaluate one-dimensional functions over equally-spaced time points ${\bf t^{\ScriptSmallRM{eval}}} = (t_1^{\ScriptSmallRM{eval}},\cdots,t_M^{\ScriptSmallRM{eval}})^T$ with step size $h$ between two consecutive points. We evaluate $\widehat G_{\ScriptSmallCap{\bm\Delta}}(s,t)$ over all possible two-dimensional grid points constructed from ${\bf t^{\ScriptSmallRM{eval}}}$, denoted by ${\bf t^{\ScriptSmallRM{eval}}}\times{\bf t^{\ScriptSmallRM{eval}}}$. Let $\widehat G_{\ScriptSmallCap{\bm\Delta}}({\bf t^{\ScriptSmallRM{eval}}}\times{\bf t^{\ScriptSmallRM{eval}}})$ be the evaluation matrix across all grid points. The estimates of eigenfunction $\phi_k(t)$ and $\lambda_k$ are derived as the eigenvectors and eigenvalues of $\widehat G_{\ScriptSmallCap{\bm\Delta}}({\bf t^{\ScriptSmallRM{eval}}}\times{\bf t^{\ScriptSmallRM{eval}}})$ adjusted for the step size $h$.
			
			In some cases the number of elements in $N(\ScriptSmallCap{\bm\Delta})$ is very limited. For example, if observations are collected from irregular and sparse locations, it is relatively rare for two pair of locations to have exactly the same spatial separation vector. More location pairs can be included by working with a sufficiently small neighborhood around a given $\bm\Delta$. Define $N(\bm\Delta,\bm\delta) = \{(i,j),\ s.t.\ \bm\Delta_{ij}\in B((\Delta_x,\Delta_y), \bm\delta)\ {\textrm{or}}\ \bm\Delta_{ij}\in B((-\Delta_x,-\Delta_y), \bm\delta)\}$ where $B(\bm\Delta,\bm\delta)$ is a neighborhood ``ball'' centering around $\bm\Delta$ with radius $\bm\delta$. The estimate of cross-covariance surface can be derived by replacing $N(\bm\Delta)$ with $N(\bm\Delta,\bm\delta)$ in \eqref{2d cross cov smoother}.

		\subsection{Anisotropy $\rm{\bf Mat\acute{e}rn}$ Model}
			Now we focus on estimating the parameters of the \matern model. We will first estimate the empirical correlation $\rho_k(\bm\Delta)$ for the cross- covariance surfaces and estimate the parameters of the \matern model by fitting them to the empirical correlations. Equation \eqref{spatial correlation} specifies the spatial covariance among fPC scores. Let $\lambda_k(\bm\Delta)$ be the $k$th eigenvalue of $G_{\ScriptSmallCap{\bm\Delta}}(s,t)$. For covariance surface, we use $\lambda_k((0,0))$ and $\lambda_k$ interchangeably. Then we have
			\begin{equation}
				G_{\ScriptSmallCap{\bm\Delta}}(s,t) = \bm\phi(s)^T{\textrm{diag}}\left(\rho^*_1(\bm\Delta)\lambda_1, \rho^*_2(\bm\Delta)\lambda_2, \cdots, \rho^*_K(\bm\Delta)\lambda_K\right)\bm\phi(t).\\
			\end{equation}
			If we further assume $\rho^*_1(\bm\Delta)\lambda_1 > \rho^*_2(\bm\Delta)\lambda_2 > \cdots > \rho^*_K(\bm\Delta)\lambda_K > 0$, then the sequence $\left\{\rho^*_k(\bm\Delta)\lambda_k\right\}_{k=1}^K$ are eigenvalues of $G_{\ScriptSmallCap{\bm\Delta}}(s,t)$ ordered from the largest to the smallest. Note that $\rho^*_k(\bm\Delta) = 1$ if $\bm\Delta = (0,0)$ for all $k$. Thus for all $\bm\Delta > 0$, $\rho^*_k(\bm\Delta)$ can be estimated as the ratio of $k$th eigenvalues of $G_{\ScriptSmallCap{\bm\Delta}}(s,t)$ and $G_{(0,0)}(s,t)$, which can be written as
			\begin{equation}
				\hat\rho^*_k(\bm\Delta) = \dfrac{\hat\lambda_k(\bm\Delta)}{\hat{\lambda}_k},		
			\end{equation}  
			where $\hat\lambda_k(\bm\Delta)$ is the $k$th eigenvalue of $\widehat G_{\ScriptSmallCap{\bm\Delta}}$ and $\hat\lambda_k$ is the $k$th eigenvalue of $\widehat G_{(0,0)}$. 									
			Suppose empirical correlations $\left\{\hat\rho^*_k(\bm\Delta_i)\right\}_{i=1}^{m}$ are obtained for $\left\{\bm\Delta_i\right\}_{i=1}^m = \left\{\bm\Delta_1,\cdots,\bm\Delta_m\right\}$. Then
			\begin{equation}\label{input for fitting anisotropy cor}
				\left\{(\bm\Delta_1,\hat\rho^*_k(\bm\Delta_1)),(\bm\Delta_2,\hat\rho^*_k(\bm\Delta_2)),\cdots,(\bm\Delta_m,\hat\rho^*_k(\bm\Delta_m))\right\}
			\end{equation}
			are used to fit \eqref{anisotropy matern correlation} and to estimate parameters $\alpha,\delta,\zeta,\nu$. If assuming separable covariance structure, empirical correlations could be pooled across $k$ to estimate parameters of the anisotropy \matern model.
			If $\left\{B(\bm\Delta_i,\bm\delta_i)\right\}_{i=1}^m$ are used, then we select one representative vector from each $B(\bm\Delta_i,\bm\delta_i)$ as input. A sensible choice of representative vectors is just $\left\{\bm\Delta_i\right\}_{i=1}^m$, the center of $\left\{B(\bm\Delta_i,\bm\delta_i)\right\}_{i=1}^m$. When fitting \eqref{anisotropy matern correlation}, the sum of squared difference between empirical and fitted correlations over all $\bm\Delta_i$'s is minimized through numerical optimization. We adopt BFGS method in implementation. More details about the quasi-Newton method can be found in \citet{broyden1970convergence}, \citet{fletcher1970new}, \citet{goldfarb1970family} and \citet{shanno1970conditioning}.
			
\subsection{Curve Reconstruction}
		Reconstructing trajectories is an important application of the SPACE model. Curve reconstructions based on SPACE model also provides an alternative perspective of ``gap-filling'' the missing data for geoscience applications as well. Equation \eqref{approx expansion} specifies the underlying formula used to reconstruct the trajectory $X_i(t)$ for each $i$. $\{\hat{\bm\phi}_k(t)\}_{k=1}^K$ and $\hat\mu(t)$ can be derived through the process described in previous sections. The only missing element now is fPC scores $\{\xi_{ik}\}_{i=1,k=1}^{N,K}$. The best linear unbiased predictors (BLUP) \citep{henderson1950estimation} of $\xi_{ik}$ are given by 
		\begin{equation}\label{fPC score BLUP}
			\check{\xi}_{ik} = E(\xi_{ik}|\{Y_{ij}\}_{i=1,j=1}^{N,n_i}).
		\end{equation}			
		To describe the closed-form solution to equation \eqref{fPC score BLUP} and to facilitate subsequent discussions, we introduce the following notations. Write ${\bf Y}_i = (Y_i(t_{i1}),\cdots,Y_i(t_{in_i}))^T$, $\widetilde{\bf Y} = ({\bf Y}_1,\cdots,{\bf Y}_N)^T$, ${\bm\mu}_i = (\mu(t_{i1}),\cdots,\mu(t_{in_i}))^T$, $\widetilde{\bm\mu} = ({\bm\mu}_1,\cdots,{\bm\mu}_N)^T$, ${\bm\xi}_i = ({\xi_{i1}},\cdots,{\xi_{iK}})^T$, $\widetilde{\bm\xi} = ({\bm\xi}_1,\cdots,{\bm\xi}_N)^T$, ${\bm\Lambda} = diag(\lambda_1,\cdots,\lambda_K)$, $\rho^*_{ijk} = \rho^*_k(\bm\Delta_{ij})$, ${\bm\rho}_{ik} = (\rho^*_{i1k},\cdots,\rho^*_{iNk})^T$, ${\bm\rho}_k = ({\bm\rho}_{1k}^T,\cdots,{\bm\rho}_{Nk}^T)^T$, ${\bm\rho}_{ij}=diag(\rho^*_{ij1},\cdots,\rho^*_{ijK})$, $\widetilde{\bm\rho} = \left[{\bm\rho}_{ij}\right]$ where $\left[\LargerCdot{1.6}_{\ScriptSmall{ij}}\right]$ represents a matrix with $ij$th entry equal to $\LargerCdot{1.6}_{\ScriptSmall{ij}}$, and $\bm\phi_{ik} = (\phi_k(t_{i1}),\cdots,\phi_k(t_{in_i}))^T$, ${\bm\phi}_i = (\bm\phi_{i1},\cdots,\bm\phi_{iK})$, $\widetilde{\bm\phi} = diag({\bm\phi}_1,\cdots,{\bm\phi}_N)$. Note diagonalization and transpose are performed before substitution in all above notations. If assuming separable covariance, we write ${\bm\rho} = \left[\rho^*_{ij}\right]$. In addition, define $\bm\Sigma(A,B)$ as the covariance matrix between $A$ and $B$, and $\bm\Sigma(A)$ as the variance matrix of $A$. Then 
		\begin{equation}
			\bm\Sigma(\widetilde{\bm\xi},\widetilde{\bm\xi}) = \left\{\begin{array}{cc}\widetilde{\bm\rho}\LargerCdot{1.6}({\bf 1}_{{N\times N}}\otimes\Lambda)&\hspace{15pt}{\textrm{nonseparable}},\\
			{\bm\rho}\otimes\bm\Lambda&\hspace{15pt}{\textrm{separable}},\end{array}\right.
   \end{equation}
		where $\otimes$ represents Kronecker product. With Gaussian assumptions, we have
		\begin{eqnarray}\label{fPC score vector BLUP}
			\check{\widetilde{\bm\xi}} &=& E(\widetilde{\bm\xi}|\widetilde{\bf Y}) = {\bm\Sigma}(\widetilde{\bm\xi},\widetilde{\bf Y}){\bm\Sigma}(\widetilde{\bf Y},\widetilde{\bf Y})^{-1}(\widetilde{\bf Y} - \widetilde{\bm\mu})\\
			&=& \bm\Sigma(\widetilde{\bm\xi},\widetilde{\bm\xi})\widetilde{\bm\phi}^T\left(\widetilde{\bm\phi}\bm\Sigma(\widetilde{\bm\xi},\widetilde{\bm\xi})\widetilde{\bm\phi}^T+\sigma^2{\bf 1}\right)^{-1}(\widetilde{\bf Y} - \widetilde{\bm\mu})\\
			&=& \left(\sigma^2\bm\Sigma(\widetilde{\bm\xi},\widetilde{\bm\xi})^{-1}+\widetilde{\bm\phi}^T\widetilde{\bm\phi}\right)^{-1}\widetilde{\bm\phi}^T(\widetilde{\bf Y} - \widetilde{\bm\mu}),
		\end{eqnarray} 
		where the last line follows the Woodbury matrix identity. For cases where $\sum_{i=1}^N n_i > NK$, the transformation of last line suggests a way to reduce the size of matrix to be inverted. The separability assumption simplifies not only the model itself but the calculation of matrix inverse as well, noting that $\bm\Sigma(\widetilde{\bm\xi},\widetilde{\bm\xi})^{-1} = ({\bm\rho}\otimes{\bm\Lambda})^{-1} = {\bm\rho}^{-1}\otimes{\bm\Lambda}^{-1}$. By substituting all elements in \eqref{fPC score vector BLUP} with corresponding estimates, the estimate of $\widetilde{\bm\xi}$ is derived as
		\begin{equation}\label{fPC score estimation}
			\widehat{\widetilde{\bm\xi}} = \left(\begin{array}{c}\widehat{\bm\xi}_1\\\vdots\\\widehat{\bm\xi}_N\end{array}\right) 
																	 = \left\{\begin{aligned} & \left({\hat\sigma}^{\ScriptSmall{2}}\widehat{\widetilde{\bm\rho}}\LargerCdot{1.6}({\bf 1}_{{N\times N}}		
																						\otimes\widehat{\bm{\bm\Lambda}}) + 		
																						\widehat{\widetilde{\bm\phi}}^{\ScriptSmall{T}}\widehat{\widetilde{\bm\phi}}\right)^{-1}
																					  \widehat{\widetilde{\bm\phi}}^{\ScriptSmall{T}}(\widetilde{\bf Y} - \widehat{\widetilde{\bm\mu}}),\hspace{15pt} & {\textrm{nonseparable
																						}}\\
																	 &	\left({\hat\sigma}^{\ScriptSmall{2}}\widehat{\bm\rho}\otimes\widehat{\bm\Lambda} + 	
																		         \widehat{\widetilde{\bm\phi}}^{\ScriptSmall{T}}\widehat{\widetilde{\bm\phi}}\right)^{-1}
																						\widehat{\widetilde{\bm\phi}}^{\ScriptSmall{T}}(\widetilde{\bf Y} - \widehat{\widetilde{\bm\mu}}),\hspace{15pt} & {\textrm{separable}}
																						\end{aligned}
																	 \right.	
	 \end{equation} 
	The reconstructed trajectory is then given by
		\begin{equation}\label{reconstruction}
			\widehat{\bf X}_i({\bf t}^{\ScriptSmall{\textrm{eval}}}) = \widehat{\bm\mu}_i({\bf t}^{\ScriptSmall{\textrm{eval}}}) + \widehat{\bm\phi}_i({\bf t}^{\ScriptSmall{\textrm{eval}}})\widehat{\bm\xi}_i.
		\end{equation}		
		
\subsection{Model Selection}\label{subsection: model selection}
			\citet{RiceSilverman1991} proposed a leave-one-curve-out cross-validation method for data which are curves. \citet{hurvich1998smoothing} introduced a methodology for choosing smoothing parameter for any linear smoother based on an improved version of Akaike information criterion (AIC). \citet{YaoMullerWang2005} pointed out that adaptation to estimated correlations when estimating the mean function with dependent data does not lead to improvements and subjective choice are often adequate.
						
			In local linear smoothing of both cross-covariance surface and mean curve, we use the default cross-validation method employed in the $sm$ package \citep{smRpackage} of R \citep{R}. In particular, that method groups observations into bins and perform leave-one-bin-out cross-validation (LOBO). We also examine other alternatives which include a) leave-one-curve-out cross-validation of \citet{RiceSilverman1991} (LOCO), b) J-fold leave-curves-out cross-validation ($\rm{LCO_J}$) and c) the improved AIC method of \citet{hurvich1998smoothing}. In all cross-validation methods, the goodness of smoothing is assessed by squared prediction error. In simulations not reported here, LOBO method demonstrates the most consistent performance in terms of eigenvalue estimation. Sparseness and noise of observations will inflate estimated eigenvalues up whereas smoothing tends to shrink estimates down. LOBO method achieves better balance between those two competing forces. We leave a more detailed investigation of this phenomenon to future work. 
		
\begin{figure}[!h]
				\begin{minipage}[t]{0.48\linewidth}
    			\centering
    			\includegraphics[width = 0.95\linewidth]
    			{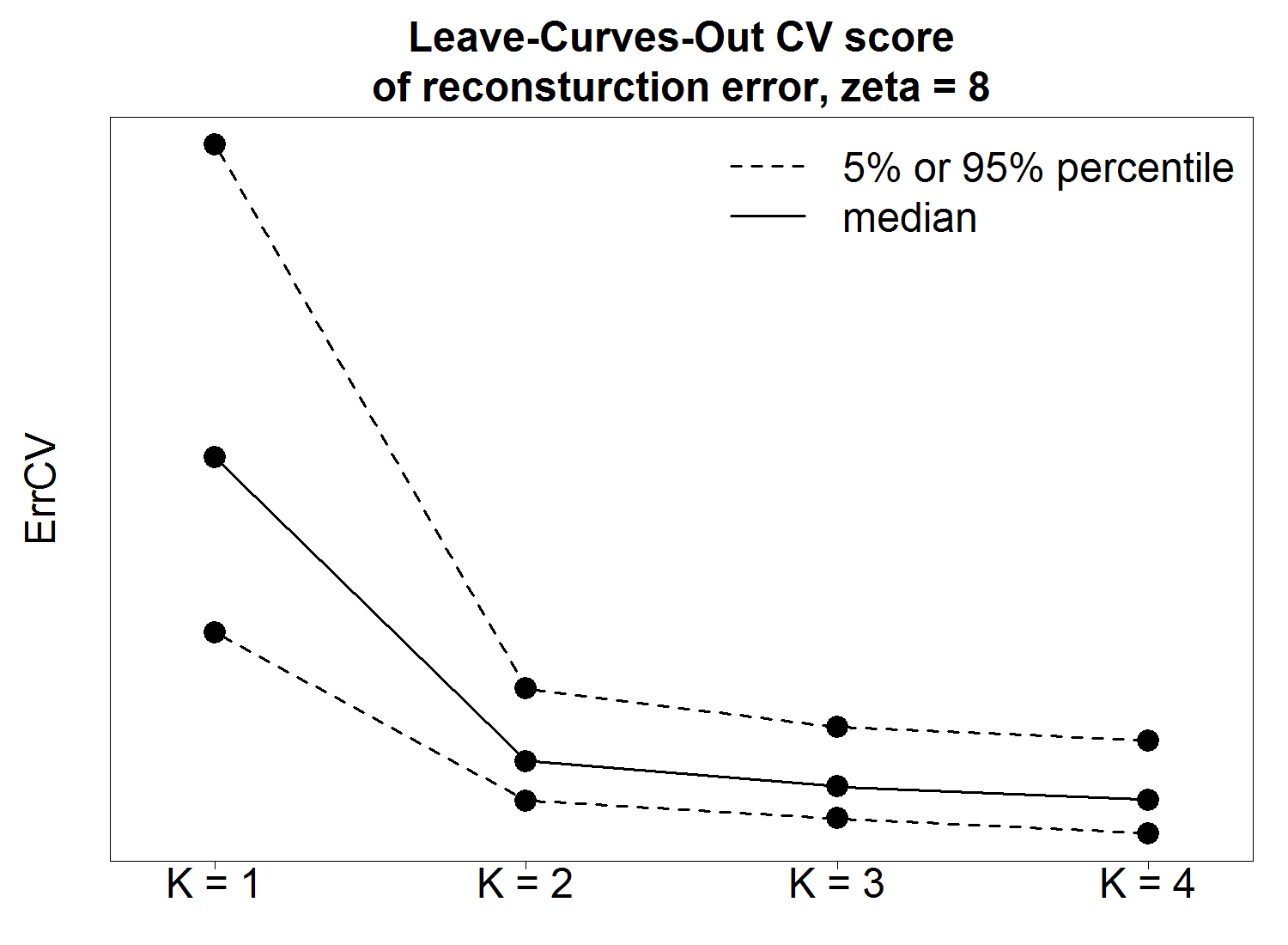}
  			\end{minipage}
  			\hspace{5pt}
				\begin{minipage}[t]{0.48\linewidth}
    			\centering
    			\includegraphics[width = 0.95\linewidth]
    			{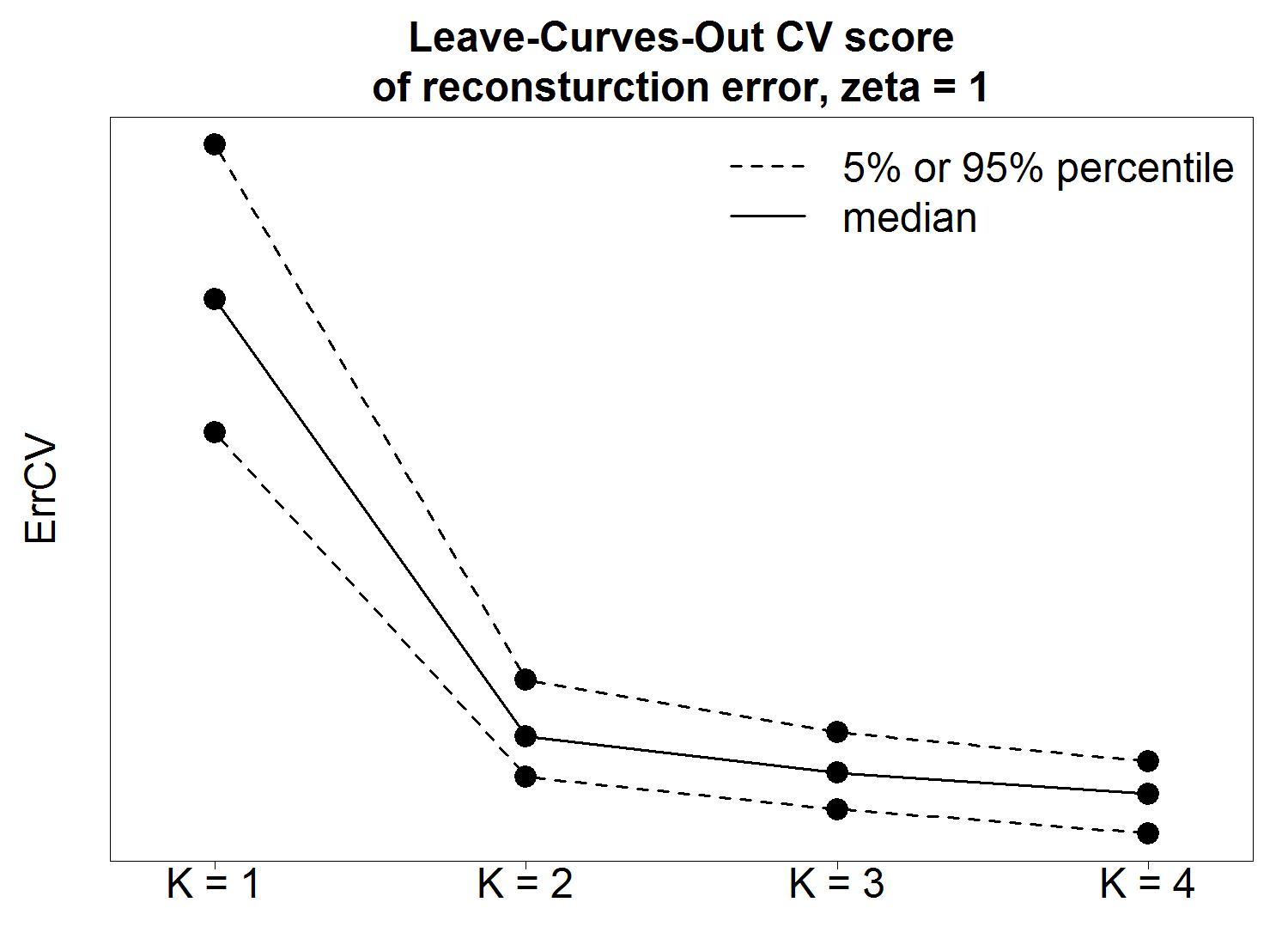}
  			\end{minipage}
 				\caption{\footnotesize{Leave-curves-out cross-validation. Simulated curves are constructed as the linear combination of eigenfunctions over time interval [0,1]. We use two eigenfunctions $\phi_1(t)\equiv 1$ and $\phi_{2}(t)=\sin(\omega t)$. Two separable spatial correlation structures are generated over a one-dimensional horizontal array of 100 equally-spaced locations. Eigenvalues are $\lambda_1 = \exp(-1)\times 10 = 3.679$ and $\lambda_2 = \exp(-2)\times 10 = 1.353$. \matern parameters are $\alpha = 0, \delta = 1, \zeta = 8, \nu = 0.5$ for high correlation scenario and $\alpha = 0, \delta = 1, \zeta = 1, \nu = 0.5$ for low correlation scenario. Corresponding correlations at distance 1 are 0.8825 and 0.3679. Noise standard deviation is $\sigma = 1$. When estimating model parameters, we use two spatial separation vectors ${\bm\Delta}_1 = (1,0)$ and ${\bm\Delta}_2 = (2,0)$. $\rm{ErrCV_5}$ is computed over 200 simulated data sets with number of eigenfunctions $K = 1,2,3,4$. Median $\rm{ErrCV_5}$ is represented by solid line and  5\% and 95\% percentiles are plotted in dashed lines.}}
  			\label{fig: KSelection}		
  		\end{figure}

 			To determine the number of eigenfunctions $K$ which sufficiently approximate the infinite dimensional process, we rely on $\rm{LCO_J}$ of curve reconstruction error which is denoted by $\rm{ErrCV_J}$ for convenience. Specifically, all curves are partitioned into $J$ complementary groups. Curves in each group serve as the testing sample. Then we select the training sample as curves that are certain distance away from testing curves to reduce the spatial correlation between the training and testing samples. For each fold, testing curves are reconstructed based on parameters estimated by corresponding training curves. Denote reconstructed curves of the $j$th fold by $\widehat{\widetilde{\bf X}}_{-j}$. Then, we compute $\rm{ErrCV_J}$ as follows,
 			\begin{equation}
 				\rm{ErrCV_J} = \sum_{j=1}^{J}\left(\widetilde{\bf X}_{-j} - \widehat{\widetilde{\bf X}}_{-j}\right)^T\left(\widetilde{\bf X}_{-j} - \widehat{\widetilde{\bf X}}_{-j}\right) \Big/ J.	
 			\end{equation}
 			\citet{Zhou2010Mixed} pointed out that cross-validation scores may keep decreasing as model complexity increases, which is also observed in our simulation studies. The 5-fold cross-validation score, $\rm{ErrCV_5}$, as a function of $K$ is illustrated in Figure \ref{fig: KSelection}. A quick drop is observed followed by a much slower decrease. Instead of finding the minimum $\rm{ErrCV_5}$, we visually select a suitable $K$ as the kink of $\rm{ErrCV_5}$ profile. In Figure \ref{fig: KSelection}, the largest drop of CV scores takes place at the correct value $K = 2$ for 100\% times out of the 200 replicated data sets.			
					
			When fitting anisotropy \matern model, $\left\{\bm\Delta_i\right\}_{i=1}^m$, the set of neighborhood structures needs to be determined. Consider the example of one-dimensional and equally-spaced locations where the $i$th location has coordinates $(i,0)$. Larger spatial separation vector corresponds to fewer pairs of locations in $N(\bm\Delta)$. Thus we often start with the most immediate neighborhood and choose $\left\{\bm\Delta_i\right\}_{i=1}^m = \left\{(1,0),(2,0),\cdots,(m,0)\right\}$ with $m$ to be decided. In simulation studies not reported here, we found that $\rm{ErrCV_J}$ is not sensitive with respect to $m$, which is especially true when spatial correlation is low. Instead of selecting an ``optimal'' $m$, we use a range of different $m$'s and take the trimmed average of estimates derived across $m$. The maximum $m$ can be selected so that the spatial correlation is too low to have meaningful impact based on prior knowledge. In simulation studies of Section \ref{section: simulation}, all estimations are made based on 20\% (each side) trimmed mean across a list of $\left\{\bm\Delta_i\right\}_{i=1}^m$'s. 
				
		\section{Asymptotic Properties of Estimates}\label{section: consistency}
		Assuming no spatial correlation, \citet{YaoMullerWang2005} demonstrated the consistency of the estimated fPC scores $\hat\xi_{ik}$, reconstructed curves $\widehat{X}_i(t)$ and other model components. Uniform convergence of the local linear estimators of mean and covariance functions on bounded intervals plays a central role in obtaining these results. In addition, the manuscript of \citet{Paul2010Manu} proposed similar nonparametric and kernel-based estimators. It was shown in their work that if the correlation between sample curves is ``weak'' in a appropriate sense, then for their estimators, the optimal rate of convergence in the correlated and i.i.d. cases are the same. Based on the framework of proof in \citet{YaoMullerWang2005} and by introducing two sufficient conditions, we are able to extend results in \citet{YaoMullerWang2005} and show the consistency of ${\hat\mu}(t)$ and cross-covariance function ${\widehat G}_{\ScriptSmall{\bm\Delta}}(s,t)$ in the following theorem.	
		
		{\bf Theorem} ${\bf 1^*}$ Under (A1.1)-(A4) and (B1.1)-(B2.2b) in \citet{YaoMullerWang2005} with $\nu=0,l=2$ in (B2.2a), $\nu=(0,0),l=2$ in (B2.2b), and (D1)-(D2) that are described in Appendix \ref{appendix: consistency},  
		$$
			\sup_{t,s\in{\cal T}}|\widehat{G}_{\ScriptSmall{\bm\Delta}}(s,t) - G_{\ScriptSmall{\bm\Delta}}(s,t)| = O_p\left(\dfrac{1}{\sqrt{|N(\bm\Delta)|}h_G^2}\right).
		$$
						
		The proof of Theorem 1* and more discussions can be found in Appendix \ref{appendix: consistency}, \ref{appendix: proofs} and \ref{appendix: intervals and regions}. In appendices, we first review the main theorems and their proofs discussed in \citet{YaoMullerWang2005}. Then, we highlight necessary modifications to accommodate the spatial dependence. 

  \section{Separability and Isotropy Tests}\label{section: separability and isotropy tests}
  	Separability of spatial covariance is a convenient assumption under which $Cov({\bm\xi}_i, {\bm\xi}_j)$ is simply a rescaled identity matrix and thus a parsimonious model could be fitted. However, we would like to design a hypothesis test to examine if this assumption is valid and justified by the data. Spatial correlation could depend not only on distance but angle as well. Whether correlation is isotropic or not may be an interesting question for researchers and is informative for subsequent analysis. In this section, we propose two hypothesis tests to address these issues based on SPACE model. 
  	\subsection{Separability Test}\label{subsection: separability test}
  		Recall the correlation matrix of the $k$th fPC scores among curves is denoted by ${\bm\rho}_k$. We model ${\bm\rho}_k$ through anisotropy \matern model through parameters $\alpha_k,\delta_k,\zeta_k$ and $\nu_k$. Suppose we use $K$ eigenfunction to approximate the underlying process. Then we partition the set $\{1,2,\cdots,K\}$ into $A$ mutually exclusive and collectively exhaustive groups ${\bf K} = \{{\cal K}_1,\cdots, {\cal K}_A\}$. Denote the number of $k$'s in group ${\cal K}_i$ by $|{\cal K}_i|$. A generic hypothesis can be formulated in which parameters across $k$'s are assumed to be the same within each group but can be different between groups. Consider a simple example where we believe correlation structures are the same among the first three fPC scores  versus the hypothesis that correlations are all different. Then the partition associated with the null hypothesis is ${\bf K}^0:\{{\cal K}_1=\{k_1,k_2,k_3\},{\cal K}_2=\{k_4\},\cdots,{\cal K}_{K-2}=\{K\}\}$. The partition for alternative hypothesis is ${\bf K}^1:\{{\cal K}_1=\{k_1\},{\cal K}_2=\{k_2\},\cdots,{\cal K}_{K}=\{K\}\}$. Both the null and alternative can take various forms of partitions. However, for illustration purpose, we consider the following test where no constraints are imposed in the alternative:
  		\begin{align}
  			& {\bf H}_0: \alpha_{i_r}=\alpha_{j_r},\delta_{i_r}=\delta_{j_r},\zeta_{i_r}=\zeta_{j_r},\nu_{i_r}=\nu_{j_r},\forall i_r,j_r\in{\cal K}_r,\ \forall r\ s.t.\ |{\cal K}_r| > 1\label{separability null}\\
  			& {\bf H}_1: \alpha_k,\delta_k,\zeta_k,\nu_k\ {\rm are\ arbitrary}\label{separability alternative}
  		\end{align}
 			The key step of the test is to construct hypothesized null curves from observed data. If sample curves are independent, hypothesized null curves can be constructed through bootstrapping fPC scores, see \citet{liu2012functional} and \citet{Li2011}. With correlated curves, reshuffling of the fPC scores is not appropriate as correlation structure would break down. However, correlated fPC scores can be transformed to uncorrelated z-scores based on the covariances estimated from anisotropy \matern model. Then z-scores are reshuffled by curve-wise bootstrapping. Next, reshuffled z-scores are transformed back to the original fPC score space based on the covariances estimated under the null hypothesis. Then hypothesized null curves are constructed as the linear combination of eigenfunctions weighted by resampled fPC scores. Let ${\bm\theta}_k\triangleq\{\alpha_k,\delta_k,\zeta_k,\nu_k\}$ be the set of anisotropy \matern parameters. The detailed procedure is described as follows.
 			\begin{mydescription}{38pt}
 				\item[{\it step 1.}] Estimate model parameters assuming they are arbitrary. Denote the estimates by ${\hat\sigma}^2,\widehat{\bm\phi}(t),\widehat{\bm\lambda}$ and $\widehat{\bm\theta}_k$. Then compute the observed test statistics ${\cal S}$ using these estimates. 
 				\item[{\it step 2.}] Estimate fPC scores using ${\hat\sigma}^2,\widehat{\bm\phi}(t),\widehat{\bm\lambda}$ and $\widehat{\bm\theta}_k$ assuming nonseparable case in \eqref{fPC score estimation}. Note $\widehat{\widetilde{\bm\xi}}$ in \eqref{fPC score estimation} is arranged by curve index $i$. We rearrange it by fPC index $k$ and let $\widehat{\bm\xi}_k$ be the vector of estimated $k$th fPC scores across all locations.    
 				\item[{\it step 3.}] Estimate ${\bm\theta}_k$ assuming ${\bf H}_0$. Specifically, for each $r$, fit anisotropy \matern model with pooled inputs $\{\bm\Delta_l,\hat\rho_k^*({\bm\Delta}_l),\forall k\in{\cal K}_r\}_{l=1}^L$.	Let $\widehat{\bm\theta}_k^{\hspace{0.5pt}\ScriptSmall{0}}$ be the pooled estimates. Note $\widehat{\bm\theta}_{k_1}^{\hspace{0.5pt}\ScriptSmall{0}} = \widehat{\bm\theta}_{k_2}^{\hspace{0.5pt}\ScriptSmall{0}}$ if $k_1,k_2$ are both in ${\cal K}_r$. 
 				\item[{\it step 4.}] For each $k$ whose associated ${\cal K}_r$ has more than one indices, let $\widehat{{\bm\rho}}_k$ be the spatial correlation matrix constructed using $\widehat{\bm\theta}_k^{\hspace{0.5pt}\ScriptSmall{0}}$. Suppose $\widehat{{\bm\rho}}_k$ has eigen-decomposition ${\bf VCV^{\ScriptSmall{T}}}$. Then $\widehat{Cov}({\bm\xi}_k) = \hat{\lambda}_k{\bf VCV^{\ScriptSmall{T}}}$. Define transformation matrix ${\bf T} = ({\hat\lambda}_k{\bf C})^{-1/2}{\bf V}^{\ScriptSmall{T}}$ and the transformed scores ${\widehat{\bf z}}_k = {\bf T}\widehat{\bm\xi}_k$. Then $\widehat{Cov}({\widehat{\bf z}}_k) = {\bf I}_{N\times N}$.
 				\item[{\it step 5.}] Resample ${\widehat{\bf z}}_k$'s through curve-wise bootstrapping. Specifically, randomly sample the curve indices $\{1,\cdots,N\}$ with replacement and let $P(i)$ be the permutated index for curve $i$. Then the $b$th bootstrapped score for curve $i$ is obtained as ${\hat{z}}^{\ScriptSmall{b}}_{ik} = {\hat{z}}_{{\ScriptSmall{P(i)}}k}$ for all $k$'s such that $|{\cal K}_r| > 1$.
 				\item[{\it step 6.}] Transform ${\widehat{\bf z}}^{\ScriptSmall{b}}_k$'s back to the space of fPC scores. Define resampled fPC scores ${\widehat{\bm\xi}}^{\ScriptSmall{b}}_k = {\bf T}^{-1}{\widehat{\bf z}}^{\ScriptSmall{b}}_k$. 
				\item[{\it step 7.}] Generate the $b$th set of Gaussian random noises $\{\epsilon^{\ScriptSmall{b}}(t_{ij})\}_{i=1,j=1}^{N,n_i}$ based on the noise standard deviation ${\hat\sigma}^2$ estimated in step 1.
 				\item[{\it step 8.}] Construct resampled observations as below
 					\begin{align}
 						Y_{i}^{\ScriptSmall{b}}(t_{ij}) & = \hat{\mu}(t_{ij}) + \sum_{|{\cal K}_r| > 1}\sum_{k\in{\cal K}_r}\hat{\xi}^{\ScriptSmall{b}}_{ik}\hat{\phi}_k(t_{ij}) + \sum_{|{\cal K}_r| = 1}\sum_{k\in{\cal K}_r}\hat{\xi}_{ik}\hat{\phi}_k(t_{ij}) + \epsilon^{\ScriptSmall{b}}(t_{ij})\label{bootstrapped curves}
 					\end{align}
 				\item[{\it step 9.}] Given $\{Y_{i}^{\ScriptSmall{b}}(t_{ij})\}_{i=1,j=1}^{N,n_i}$, estimate anisotropy \matern parameters. Denote the estimates by $\widehat{\bm\theta}^{\ScriptSmall{b}}_k$. Then compute test statistics ${\cal S}^{\ScriptSmall{b}}$.
 				\item[{\it step 10.}] Repeat steps 5 to 9 for B times and obtain $\{{\cal S}^{\ScriptSmall{b}}\}_{b=1}^{B}$ which form the empirical null distribution of the test statistics. Then make decision by comparing the null distribution with the observed test statistics ${\cal S}$.
			\end{mydescription}			
 			
 			 Note that an alternative and theoretically more precise way of doing step 4 is to use covariances constructed based on $\widehat{\bm\theta}_k$ as opposed to $\widehat{\bm\theta}_k^{\hspace{0.5pt}\ScriptSmall{0}}$. The alternative way is assumed to remove spatial correlation more thoroughly but simulation results not reported suggest that it introduces more volatility into z-scores that offsets the benefit of lower residual spatial correlation in them. For testing the ${\bf H}_0$ above, we transform $\widehat{\bm\theta}_k$ and $\widehat{\bm\theta}^{\ScriptSmall{b}}_k$ from parameter space to correlation space. Let $\bar{\rho}^*_r({\bm\Delta}^{\ScriptSmallRM{eval}}; \widehat{\bm\theta}_k) = \sum_{k\in{\cal K}_r}\rho^*({\bm\Delta}^{\ScriptSmallRM{eval}}; \widehat{\bm\theta}_k)/{|{\cal K}_r|}$ be the average correlation computed at separation vector ${\bm\Delta}^{\ScriptSmallRM{eval}}$ across $k$'s in ${\cal K}_r$. Then define test statistics as 
 			\begin{eqnarray}
 				{\cal S} &=& \sum_{r=1}^A\left(\dfrac{\sum_{k\in{\cal K}_r}\left(\rho^*({\bm\Delta}^{\ScriptSmallRM{eval}}; \widehat{\bm\theta}_k) - \bar{\rho}^*_r({\bm\Delta}^{\ScriptSmallRM{eval}}; \widehat{\bm\theta}_k)\right)^2}{|{\cal K}_r|}\right)^{1/2}\\
				{\cal S}^{\ScriptSmall{b}} &=& \sum_{r=1}^A\left(\dfrac{\sum_{k\in{\cal K}_r}\left(\rho^*({\bm\Delta}^{\ScriptSmallRM{eval}}; \widehat{\bm\theta}^{\ScriptSmall{b}}_k) - \bar{\rho}^*_r({\bm\Delta}^{\ScriptSmallRM{eval}}; \widehat{\bm\theta}^{\ScriptSmall{b}}_k)\right)^2}{|{\cal K}_r|}\right)^{1/2}.	
 			\end{eqnarray}
			We tried other statistics and found that statistics in the correlation space works generally better than that in the parameter space. In our implementations, we choose ${\bm\Delta}^{\ScriptSmallRM{eval}}$ to be the most immediate neighborhood in Euclidean distance. 
			
			\subsection{Isotropy Test}\label{subsection: isotropy test}
			Isotropy test is essentially the same as the separability test. To make the test more general, we assume a prior that correlation structures are the same across $k$'s within ${\cal K}_r$ for $r = 1,2,\cdots,A$. Note that within each ${\cal K}_r$, $\alpha_k$'s and $\delta_k$'s must take equal values across $k$'s. Let ${\bf K}_a = \{{\cal K}_{r_1},\cdots,{\cal K}_{r_a}\}$ be the set whose $\alpha_k$'s are all assumed to be zeros. Then let ${\bf K}_{a}^c$ be the complement set of ${\bf K}_{a}$ where $1 \leq a\leq A$. To test the following hypothesis
			\begin{align}
  			& {\bf H}_0: \exists{\bf K}_a = \{{\cal K}_{r_1},\cdots,{\cal K}_{r_a}\}\ s.t.\ \alpha_k = 0,\ \forall k\in{\cal K}_{r_i},\ i = 1,\cdots,a\label{isotropy null}\\
  			& {\bf H}_1: \forall k,\ \alpha_k\neq 0,\label{isotropy alternative}
  		\end{align} 	  
  		we propose a similar procedure which slightly differs from the separability test in step 1, step 3, step 7 and step 8. Specifically, we replace those steps in separability tests with the following. 
 			\begin{mydescription}{44pt}
 				\item[{\it step 1*.}] Estimate model parameters under the prior. In particular, $\widehat{\bm\theta}_k$ is obtained by pooling empirical correlations from all $k$'s within ${\cal K}_r$. Then compute ${\cal S}$ based on these estimates.
 				\item[{\it step 3*.}] Under the null hypothesis, we estimate anisotropy \matern parameters. For any $k\in{\cal K}_r$ that belongs to ${\bf K}_a$, $\hat\alpha_k^{\ScriptSmall{0}}$ and $\hat\delta_k^{\ScriptSmall{0}}$ are fixed at zero and one respectively whereas ${\hat\zeta}_k^{\ScriptSmall{0}}$ and ${\hat\nu}_k^{\ScriptSmall{0}}$ are estimated from $\{\bm\Delta_l,\hat\rho_k^*({\bm\Delta}_l),\forall k\in{\cal K}_r\}_{l=1}^L$. Then $\widehat{\bm\theta}_k^{\hspace{0.5pt}\ScriptSmall{0}} = (0,1,{\hat\zeta}_k^{\ScriptSmall{0}},{\hat\nu}_k^{\ScriptSmall{0}})$ for any $k\in{\cal K}_r$ that belongs to ${\bf K}_a$. For any $k\in{\cal K}_r$ that belongs to ${\bf K}_a^c$, no extra estimation is needed. We keep fPC scores estimated in step 2 for those $k$'s and perform no transformation on them.
 			  \item[{\it step 8*.}] Construct resampled curves as below
 					\begin{align}
 						Y_{i}^{\ScriptSmall{b}}(t_{ij}) & = \hat{\mu}(t_{ij}) + \sum_{{\cal K}_r\in{\bf K}_a}\sum_{k\in{\cal K}_r}\hat{\xi}^{\ScriptSmall{b}}_{ik}\hat{\phi}_k(t_{ij}) + \sum_{{\cal K}_r\in{\bf K}^c_a}\sum_{k\in{\cal K}_r}\hat{\xi}_{ik}\hat{\phi}_k(t_{ij}) + \epsilon^{\ScriptSmall{b}}(t_{ij})\label{bootstrapped curves isotropy}
 					\end{align}
 				\item[{\it step 9*.}] Given $\{Y_{i}^{\ScriptSmall{b}}(t_{ij})\}_{i=1,j=1}^{N,n_i}$, estimate anisotropy \matern parameters .
			\end{mydescription}
			Accordingly, the test statistics ${\cal S}$ could be as simple as the sum of absolute $\hat\alpha$'s, 
			\begin{equation}
				{\cal S} = \sum_{{\cal K}_r\in{\bf K}_a}|{\hat\alpha}_r|	
			\end{equation} 
			where ${\hat\alpha}_r$ represents the common estimate for all $k$'s in set ${\cal K}_r$. Testing procedures described above can be easily extended to the case where both separability and isotropy effect are of interest. 
  
 \section{Simulations}\label{section: simulation} 
 	In this section we present simulation studies of the SPACE model. We examine the estimation of anisotropy \matern parameters and curve reconstructions in Subsection \ref{subsection: sims model estimation}. Simulation results of hypothesis tests are shown in Subsetction \ref{subsection: sims hypothesis tests}. For both one-dimensional and two-dimensional locations, we consider grid points with integer coordinates. Specifically, two-dimensional coordinates take the form of $\{(i,j), 1\leq i,j\leq E\}$ where $E$ is referred to as the edge length. Similarly, one-dimensional locations are represented by $\{(0,1),(0,2),\cdots,(0,N)\}$ where $N$ is the total number of locations. For simplicity, we examine a subset of anisotropy \matern family by fixing $\nu$ at 0.5, which corresponds to the spatial autoregressive model of order one. When constructing simulated curves, we employ two eigenfunctions which are $\phi_1(t) = 1$ and $\phi_2(t) = \sin(2\pi t)$, with eigenvalues $\lambda_1 = \exp(-1)\times 10 \approx 3.68$ and $\lambda_2 = \exp(-2)\times 10\approx 1.35$, and the mean function is $\mu(t) = 0$. All functions are built on the closed interval $[0,1]$. On each curve, 10 observations are generated at time points randomly selected from 101 equally-spaced points on [0,1]. The spatial correlation between fPC scores are generated by anisotropy \matern model in equation \eqref{anisotropy matern correlation}. In addition, as indicated in Subsection \ref{subsection: model selection}, we estimate \matern parameters over a list of nested $\left\{\bm\Delta_i\right\}_{i=1}^m$'s. In the one-dimensional case, we examine the following list: $\{\{\bm\Delta_i\}_{i=1}^m\}_{m=1}^{20} = \{(0,1),\cdots,(0,m)\}_{m=1}^{20}$. In the two-dimensional case, the entire set of $\bm\Delta$ used for estimation is defined as $\{\bm\Delta_i\}_{i=1}^{24} =$ \{(1,0), (1,1), (0,1), (1,-1), (2,0), (2,1), (2,2), (1,2), (0,2), (1,-2), (2,-2), (2,-1), (3,0), (3,1), (3,2), (3,3), (2,3), (1,3), (0,3), (1,-3), (2,-3), (3,-3), (3,-2), (3,-1)\}. Estimation is performed over $\{\bm\Delta_1,\cdots,\bm\Delta_{m+4}\}_{m=1}^{20}$. Final estimates are derived as the 20\% (each side) trimmed mean across $m$ and the order of $\bm\Delta_i$ doesn't have meaningful impact on the estimates.
	
	\subsection{Model Estimation}\label{subsection: sims model estimation}
 			We examine the estimation of SPACE model in both one-dimensional and two-dimensional cases. In one-dimensional case, we are interested in the estimation of $\zeta$ whereas the estimation of $\alpha$	is our main focus in two-dimensional case. Results are summarized in Table \ref{table: 1DEst} and \ref{table: 2DEst}. Note the first order derivative of \matern correlation with respect to $\zeta$ is flattened as the correlation approaches 1. Thus more extreme large estimates of $\zeta$ are expected, which leads to the positive skewness observed in Table \ref{table: 1DEst}. We also look at the estimation performance in the correlation space. Specifically, we examine the estimated correlation at $\bm\Delta^{\ScriptSmallRM{eval}} = (0,1)$. Estimation is more difficult in two-dimensional case where two more parameters $\alpha$ and $\delta$ need to be estimated. In general, estimates based on more significant fPCs have better quality in terms of root-mean-squared error (RMSE). To assess the curve reconstruction performance, consider a grid of time points for function evaluation ${\bf t^{\ScriptSmallRM{eval}}} = (t^{\ScriptSmallRM{eval}}_1,t^{\ScriptSmallRM{eval}}_2,\cdots,t^{\ScriptSmallRM{eval}}_n)$. Let ${\rm{Err}} = \dfrac{1}{nN}\sum_{i=1}^N\sum_{j=1}^n(X_i(t^{\ScriptSmallRM{eval}}_j) - {\hat X}_i(t^{\ScriptSmallRM{eval}}_j))^2$ be the curve reconstruction error. Then define the reconstruction improvement (IP) as ${\rm IP} = \log({p\rm{Err}}_{\ScriptSmallRM{PACE}}\big/{\rm{Err}}_{\ScriptSmallRM{SPACE}})$. Out of 100 simulated data sets, we calculate the percentage of IP greater than 0. The improvement of SPACE over PACE is more prominent in scenarios with larger noise and higher spatial correlation. It is easy to show that when $\sigma = 0$ and the number of eigenfunctions is greater than the maximum number of observations per curve, SPACE produces exactly the same reconstructed curves as PACE. If noise is large, information provided by each curve itself is more contaminated by noise relative to neighboring locations which provide more useful guidance in curve reconstruction. As spatial correlation increases, observations at neighboring locations are more informative.

\setlength{\tabcolsep}{.5\tabcolsep}

		\begin{table}[ht]
	\caption{\small{Model parameter estimation for one-dimensional locations. Summary statistics of estimated $\zeta$ and $\rho^*(\bm\Delta^{\ScriptSmallRM{eval}})$  from the first and second fPCs in 6 scenarios are presented.}}				
				\label{table: 1DEst}
				\begin{center}
				\scalebox{0.9}{
				\begin{tabular}{cccccccccccccc}
					\hline
					\hline
					\multicolumn{1}{c|}{\multirow{2}{*}{Scenario}} & \multicolumn{4}{c|}{Setting}  & \multicolumn{4}{c|}{$\hat\zeta$} & \multicolumn{4}{c|}{$\hat\rho$} & \multirow{2}{*}{${\rm \% (IP > 0)}$}\\
					\cline{2-13}
		      \multicolumn{1}{c|}{}			                     & $\sigma$ & fPC & $\zeta$ & \multicolumn{1}{c|}{$\rho$} & mean & median & std & \multicolumn{1}{c|}{RMSE} & mean & median & std & \multicolumn{1}{c|}{RMSE} & \\\hline
					\multicolumn{1}{c|}{\multirow{2}{*}{separable 1}} & \multirow{2}{*}{0.2} & 1st & 5 & \multicolumn{1}{c|}{0.82} & 5.43 & 4.96 & 1.90 & \multicolumn{1}{c|}{1.92} & 0.818 & 0.821 & 0.049 & \multicolumn{1}{c|}{0.050} & \multirow{2}{*}{63\%}\\
					\multicolumn{1}{c|}{}														  &								       & 2nd & 5 & \multicolumn{1}{c|}{0.82} & 4.38 & 4.14 & 1.55 & \multicolumn{1}{c|}{1.66} & 0.778 & 0.785 & 0.059 & \multicolumn{1}{c|}{0.072} & \\\hline
					\multicolumn{1}{c|}{\multirow{2}{*}{separable 2}} & \multirow{2}{*}{1}   & 1st & 5 & \multicolumn{1}{c|}{0.82} & 5.33 & 4.98 & 1.92 & \multicolumn{1}{c|}{1.94} & 0.813 & 0.816 & 0.051 & \multicolumn{1}{c|}{0.052} & \multirow{2}{*}{99\%}\\
					\multicolumn{1}{c|}{} 														&								     	 & 2nd & 5 & \multicolumn{1}{c|}{0.82} & 4.20 & 3.87 & 1.32 & \multicolumn{1}{c|}{1.54} & 0.774 & 0.772 & 0.053 & \multicolumn{1}{c|}{0.069} & \\\hline
					\multicolumn{1}{c|}{\multirow{2}{*}{separable 3}} & \multirow{2}{*}{0.2} & 1st & 2 & \multicolumn{1}{c|}{0.61} & 2.46 & 2.34 & 0.73 & \multicolumn{1}{c|}{0.86} & 0.648 & 0.651 & 0.079 & \multicolumn{1}{c|}{0.089} & \multirow{2}{*}{58\%}\\
					\multicolumn{1}{c|}{}														  &								       & 2nd & 2 & \multicolumn{1}{c|}{0.61} & 2.53 & 2.40 & 0.75 & \multicolumn{1}{c|}{0.91} & 0.654 & 0.659 & 0.084 & \multicolumn{1}{c|}{0.096} & \\\hline
					\multicolumn{1}{c|}{\multirow{2}{*}{separable 4}} & \multirow{2}{*}{1}   & 1st & 2 & \multicolumn{1}{c|}{0.61} & 2.43 & 2.33 & 0.77 & \multicolumn{1}{c|}{0.87} & 0.642 & 0.650 & 0.084 & \multicolumn{1}{c|}{0.091} & \multirow{2}{*}{97\%}\\
					\multicolumn{1}{c|}{}														  &								     	 & 2nd & 2 & \multicolumn{1}{c|}{0.61} & 2.37 & 2.19 & 0.71 & \multicolumn{1}{c|}{0.80} & 0.637 & 0.634 & 0.084 & \multicolumn{1}{c|}{0.089} & \\\hline
					\multicolumn{1}{c|}{\multirow{2}{*}{non-separable 1}} & \multirow{2}{*}{0.5} & 1st & 6 & \multicolumn{1}{c|}{0.85} & 5.04 & 4.53 & 2.52 & \multicolumn{1}{c|}{2.67} & 0.798 & 0.818 & 0.075 & \multicolumn{1}{c|}{0.092} & \multirow{2}{*}{74\%}\\					
					\multicolumn{1}{c|}{}														  &								       & 2nd & 2 & \multicolumn{1}{c|}{0.61} & 1.54 & 1.55 & 0.51 & \multicolumn{1}{c|}{0.68} & 0.504 & 0.513 & 0.116 & \multicolumn{1}{c|}{0.151} & \\\hline					
					\multicolumn{1}{c|}{\multirow{2}{*}{non-separable 2}} & \multirow{2}{*}{1} & 1st & 6 & \multicolumn{1}{c|}{0.85} & 5.07 & 4.50 & 2.55 & \multicolumn{1}{c|}{2.71} & 0.793 & 0.804 & 0.083 & \multicolumn{1}{c|}{0.102} & \multirow{2}{*}{100\%}\\					
					\multicolumn{1}{c|}{}														  &								       & 2nd & 2 & \multicolumn{1}{c|}{0.61} & 1.48 & 1.50 & 0.46 & \multicolumn{1}{c|}{0.65} & 0.505 & 0.511 & 0.102 & \multicolumn{1}{c|}{0.143} & \\\hline					
					\hline					
				\end{tabular}
				}
				\end{center}
				\vspace{0pt}
			
			\end{table}				
  		
			\begin{table}[ht]

	\caption{\small{Model parameter estimation for two-dimensional locations. Summary statistics of estimated $\alpha$ and $\rho^*(\bm\Delta^{\ScriptSmallRM{eval}})$  from the first and second fPCs in 5 scenarios are presented. $\sigma = 1$ in all 5 scenarios.}}				
				\label{table: 2DEst}				\begin{center}

				\scalebox{0.9}{
				\begin{tabular}{cccccccccccccc}
					\hline
					\hline
					\multicolumn{1}{c|}{\multirow{2}{*}{Scenario}} & \multicolumn{4}{c|}{Setting}  & \multicolumn{4}{c|}{$\hat\alpha$} & \multicolumn{4}{c|}{$\hat\rho$} & \multirow{2}{*}{${\rm \% (IP > 0)}$}\\
					\cline{2-13}
		      \multicolumn{1}{c|}{}			                     & fPC & $\zeta$ & $\alpha$ & \multicolumn{1}{c|}{$\rho$} & mean & median & std & \multicolumn{1}{c|}{RMSE} & mean & median & std & \multicolumn{1}{c|}{RMSE} & \\\hline
					\multicolumn{1}{c|}{\multirow{2}{*}{separable 1}} & 1st & 6 & 30 & \multicolumn{1}{c|}{0.66} & 28.35 & 28.75 & 4.42 & \multicolumn{1}{c|}{4.70} & 0.505 & 0.522 & 0.111 & \multicolumn{1}{c|}{0.194} & \multirow{2}{*}{85\%}\\
					\multicolumn{1}{c|}{}														  & 2nd & 6 & 30 & \multicolumn{1}{c|}{0.66} & 29.26 & 29.54 & 6.56 & \multicolumn{1}{c|}{6.57} & 0.425 & 0.426 & 0.111 & \multicolumn{1}{c|}{0.263} & \\\hline
					\multicolumn{1}{c|}{\multirow{2}{*}{separable 2}} & 1st & 6 & 60 & \multicolumn{1}{c|}{0.79} & 61.19 & 60.60 & 4.94 & \multicolumn{1}{c|}{5.05} & 0.668 & 0.662 & 0.095 & \multicolumn{1}{c|}{0.151} & \multirow{2}{*}{74\%}\\
					\multicolumn{1}{c|}{} 														& 2nd & 6 & 60 & \multicolumn{1}{c|}{0.79} & 59.89 & 60.31 & 9.35 & \multicolumn{1}{c|}{9.31} & 0.618 & 0.630 & 0.101 & \multicolumn{1}{c|}{0.196} & \\\hline
					\multicolumn{1}{c|}{\multirow{2}{*}{separable 3}} & 1st & 3 & 30 & \multicolumn{1}{c|}{0.44} & 30.10 & 30.62 & 5.95 & \multicolumn{1}{c|}{5.92} & 0.401 & 0.396 & 0.102 & \multicolumn{1}{c|}{0.109} & \multirow{2}{*}{65\%}\\
					\multicolumn{1}{c|}{}														  & 2nd & 3 & 30 & \multicolumn{1}{c|}{0.44} & 29.11 & 29.18 & 5.68 & \multicolumn{1}{c|}{5.72} & 0.365 & 0.363 & 0.106 & \multicolumn{1}{c|}{0.130} & \\\hline
					\multicolumn{1}{c|}{\multirow{2}{*}{separable 4}} & 1st & 3 & 60 & \multicolumn{1}{c|}{0.62} & 61.95 & 60.96 & 4.86 & \multicolumn{1}{c|}{5.21} & 0.586 & 0.584 & 0.102 & \multicolumn{1}{c|}{0.106} & \multirow{2}{*}{66\%}\\
					\multicolumn{1}{c|}{}														  & 2nd & 3 & 60 & \multicolumn{1}{c|}{0.62} & 60.23 & 60.68 & 7.32 & \multicolumn{1}{c|}{7.29} & 0.538 & 0.541 & 0.107 & \multicolumn{1}{c|}{0.133} & \\\hline
					\multicolumn{1}{c|}{\multirow{2}{*}{non-separable 1}} & 1st & 5 & 75 & \multicolumn{1}{c|}{0.85} & 66.20 & 68.84 & 10.83 & \multicolumn{1}{c|}{11.46} & 0.702 & 0.700 & 0.098 & \multicolumn{1}{c|}{0.149} & \multirow{2}{*}{77\%}\\					
					\multicolumn{1}{c|}{}														  & 2nd & 5 & 45 & \multicolumn{1}{c|}{0.61} & 51.55 & 51.88 & 11.73 & \multicolumn{1}{c|}{13.41} & 0.535 & 0.537 & 0.094 & \multicolumn{1}{c|}{0.163} & \\\hline								
					\hline					
				\end{tabular}
				}
				\end{center}
				\vspace{0pt}
			
			\end{table}	
						
 			\subsection{Hypothesis Test}\label{subsection: sims hypothesis tests}
				We evaluate the hypothesis tests proposed in Section \ref{section: separability and isotropy tests}. Separability and isotropy tests are implemented based on simulated data sets on one-dimensional and two-dimensional locations respectively. In each test, 100 curves are created in each data set and 25 data sets are generated. We first focus on separability test by examining different statistics: absolute difference in $\zeta_1$ and $\zeta_2$, difference of spatial correlation at $\bm\Delta^{\ScriptSmallRM{eval}} = (0,1)$ and absolute difference of spatial correlation at $\bm\Delta^{\ScriptSmallRM{eval}} = (0,1)$. Next, we test the isotropy effect assuming separability and the test statistics is simply $\hat\alpha$. Each alternative test is evaluated at a single set of parameters. All settings are described in Table \ref{table: Separability Test} and \ref{table: Isotropy Test}. With nominal size set to 5\%, the two-sided empirical sizes and powers are summarized in Table \ref{table: Separability Test} and \ref{table: Isotropy Test}. Both tests deliver reasonable sizes and powers. 			
			\begin{table}[ht]
	\caption{\small{Empirical sizes and powers of separability test.}}				
				\label{table: Separability Test}
				\begin{center}
				\scalebox{0.9}{
				\begin{tabular}{ccccccccc}
					\hline
					\hline
					  \multirow{2}{*}{Test} & \multicolumn{5}{c}{Settting} & \hspace{20pt} \multirow{2}{*}{$\left|\hat\zeta_1 - \hat\zeta_2\right|$} \hspace{20pt} & \multirow{2}{*}{$\left|\rho^*_1(\bm\Delta^{\ScriptSmallRM{eval}}) - \rho^*_2(\bm\Delta^{\ScriptSmallRM{eval}})\right|$} & \hspace{5pt} \multirow{2}{*}{$\rho^*_1(\bm\Delta^{\ScriptSmallRM{eval}}) - \rho^*_2(\bm\Delta^{\ScriptSmallRM{eval}})$} \hspace{5pt}\\
						\cline{2-6}
						 & fPC & $\alpha$ & $\delta$ & $\nu$ & $\zeta$ & & & \\\hline
					\multirow{2}{*}{Separability Size} & 1st & 0 & 1 & 0.5 & 6 & \multirow{2}{*}{1/25} & \multirow{2}{*}{1/25} & \multirow{2}{*}{1/25}\\
																						 & 2nd & 0 & 1 & 0.5 & 6 & & & \\				
					\multirow{2}{*}{Separability Power} & 1st & 0 & 1 & 0.5 & 8 & \multirow{2}{*}{23/25} & \multirow{2}{*}{23/25} & \multirow{2}{*}{23/25}\\
																							& 2nd & 0 & 1 & 0.5 & 1 & & & \\\hline
					\hline
				\end{tabular}
				}
				\end{center}
				\vspace{0pt}
			
			\end{table}		

			\begin{table}[ht]
		\caption{\small{Empirical sizes and powers of isotropy test.}}				
				\label{table: Isotropy Test}		
		\begin{center}
				\scalebox{0.9}{
				\begin{tabular}{c|cccccc}
					\hline
					\hline
					  \multirow{2}{*}{Test} & \multicolumn{5}{c}{Settting} & \multirow{2}{*}{$\hat\alpha$}\\
						\cline{2-6}
						 & fPC & $\alpha$ & $\delta$ & $\nu$ & $\zeta$ & \\\hline
					\multirow{2}{*}{Isotropy Size} & 1st & 0 & 1 & 0.5 & 5 & \multirow{2}{*}{1/25}\\
																						 & 2nd & 0 & 1 & 0.5 & 5 & \\				
					\multirow{2}{*}{Isotropy Power} & 1st & 30 & 8 & 0.5 & 5 & \multirow{2}{*}{22/25}\\
																							& 2nd & 30 & 8 & 0.5 & 5 &\\\hline
					\hline
				\end{tabular}
				}
				\end{center}
				\vspace{0pt}
		
			\end{table}						
		
	\section{Harvard Forest Data}\label{section: Harvard Forest}
		SPACE model is motivated by the spatial correlation observed in the Harvard Forest vegetation index data described in Section \ref{subsection: matern class} and observed in \citet{liu2012functional}. In this section, we evaluate the SPACE model and isotropy test on the Enhanced Vegetation Index (EVI) at Harvard Forest, the same data set previously examined in \citet{liu2012functional}. EVI is constructed from surface spectral reflectance measurements obtained from Moderate Resolution Imaging Spectroradiometer onboard NASA's Terra and Aqua satellites. In particular, the EVI data used in this work is extracted for a 25 \time 25 pixel window which covers approximately an area of 134 $\rm{Km}^2$. The area is centered over the Harvard Forest Long Term Experimental Research site in Petershan, MA. The data is provided at 8-day intervals over the period from 2001 to 2006. More details about the Harvard Forest EVI data can be found in \citet{liu2012functional}.   
		
		We first focus on verifying the directional effect observed in the second fPC scores through the proposed isotropy test. The Harvard Forest EVI data is observed over a dense grid of regularly spaced time points. EVI observed at each individual location is smoothed using regularization approach based on the expansion of saturated Fourier basis. Specifically, the sum of squared error loss and the penalization of total curvature is minimized. Let $x_i(t)$ be the smoothed EVI at the $i$th location and $\theta_k(t)$ be the $k$th basis function. Define the $N\times 1$ vector-valued function ${\bf x}(t) = (x_1(t),\cdots,x_N(t))^T$ and the $K\times 1$ vector-valued function $\bm\theta(t) = (\theta_1(t),\cdots,\theta_K(t))^T$. All $N$ curves can be expressed as ${\bf x}(t) = {\bf C}{\bm\theta}(t)$ where ${\bf C}$ is the coefficient matrix of size $N\times K$. The functional PCA problem reduces to the multivariate PCA of the coefficient array ${\bf C}$. Assuming $K < N$ and let ${\bf U}$ be the $K\times K$ matrix of eigenvectors of ${\bf C}^T{\bf C}$. Let ${\bm\phi}(t) = (\phi_1(t), \cdots, \phi_K(t))^T$ be the vector-valued eigenfunction and $\bm\xi$ be the $N\times K$ matrix of fPC scores. Then we have ${\bm\phi}(t) = {\bf U}^T{\bm\theta}(t)$, $\bm\xi = {\bf CU}$ and ${\bf x}(t) = {\bm\xi}{\bm\phi}(t)$. 
		
\begin{figure}[!h]
				\begin{minipage}[t]{0.48\linewidth}
    			\centering
    			\includegraphics[width = 0.95\linewidth]
    			{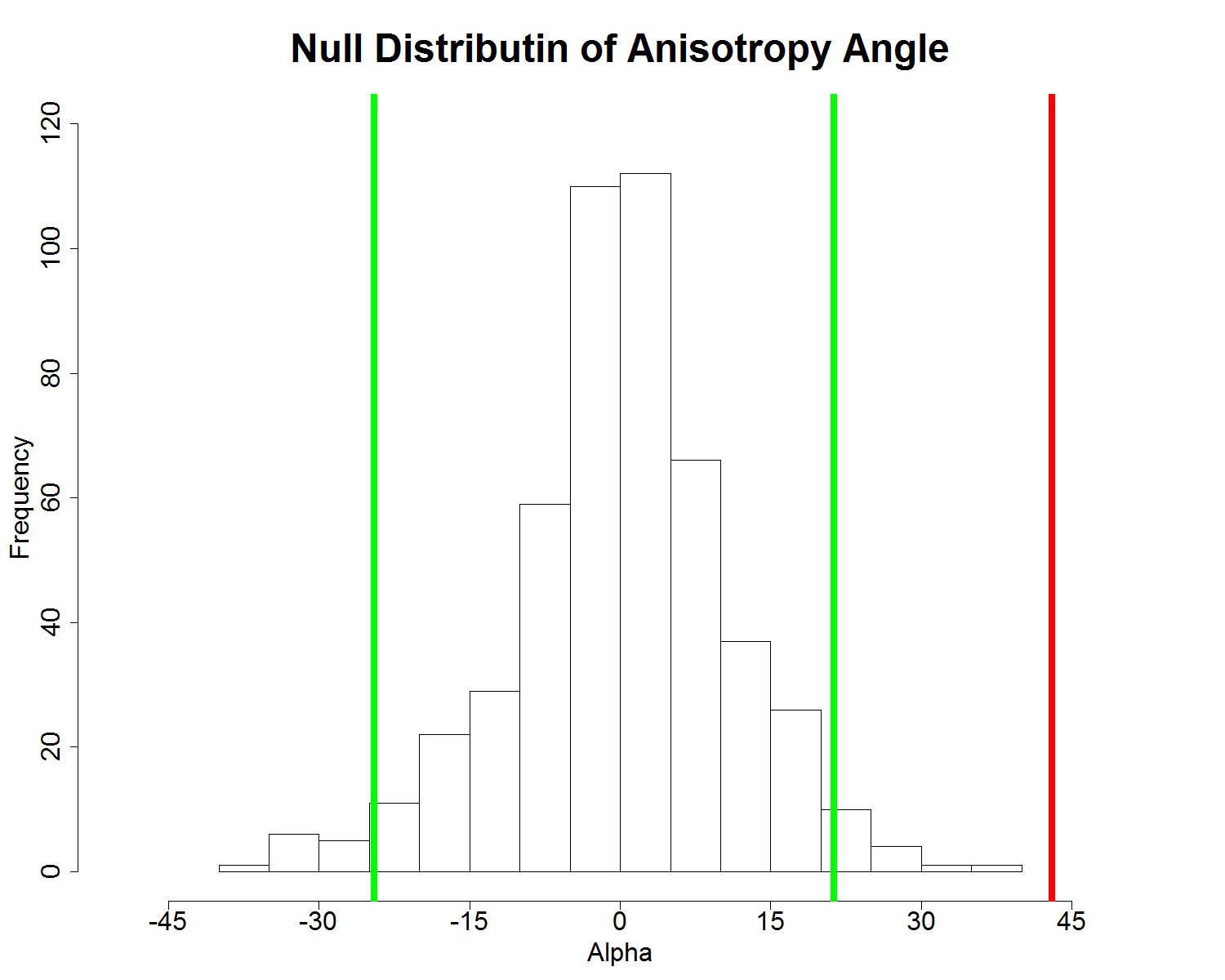}
  			\end{minipage}
  			\hspace{5pt}
				\begin{minipage}[t]{0.48\linewidth}
    			\centering
    			\includegraphics[width = 0.95\linewidth]
    			{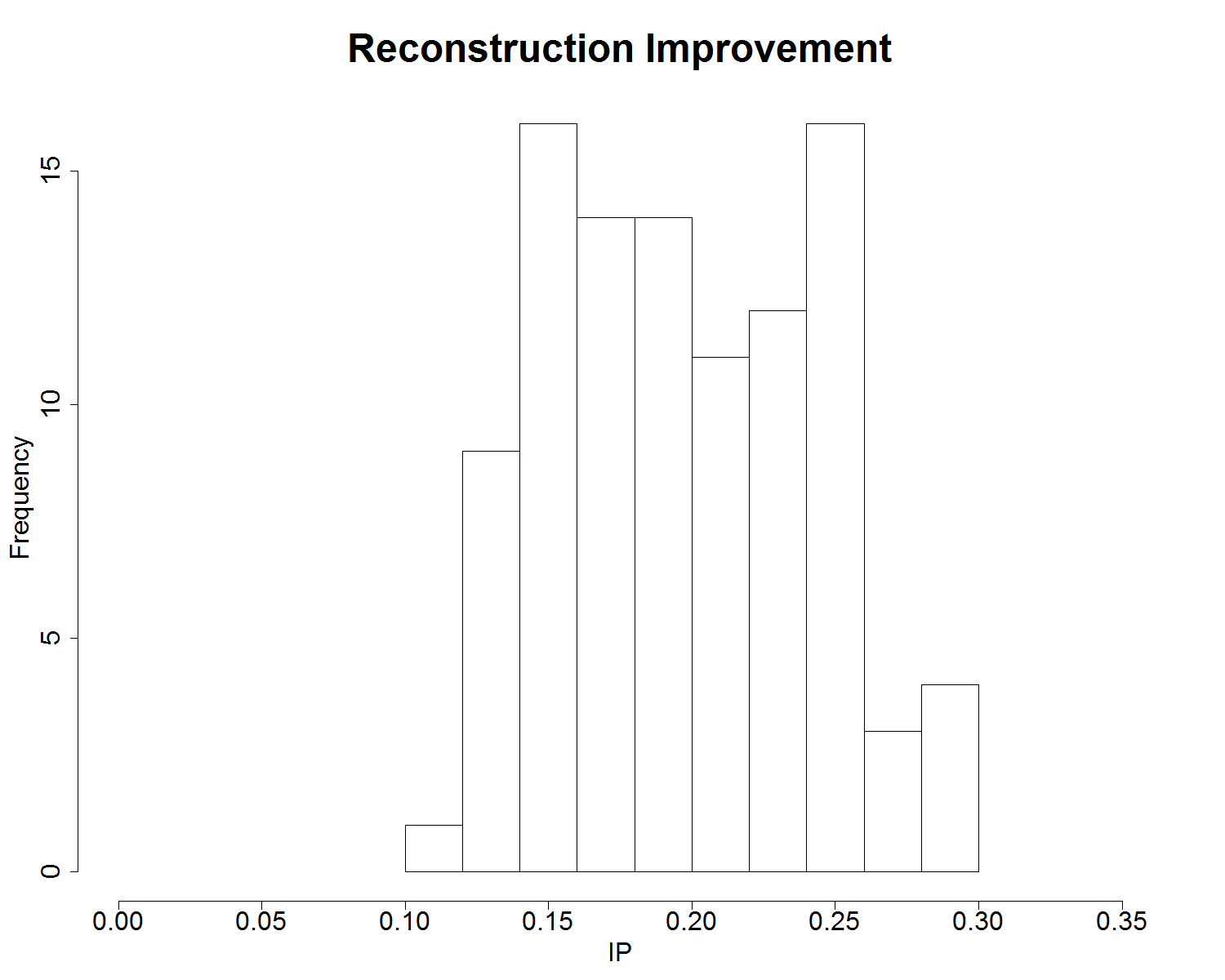}
  			\end{minipage}
 				\caption{\footnotesize{500 sets of hypothesized curves are constructed in the isotropy test. In the plot on the left, upper and lower 2.5\% percentiles of the null distribution are indicated by green vertical lines. The observed anisotropy angle is 43 degree indicated by the red vertical line. The plot on the right shows the distribution of ${\rm IP} = \log({\rm{Err}}_{\ScriptSmallRM{PACE}}\big/{\rm{Err}}_{\ScriptSmallRM{SPACE}})$. To assess the gap-filling performance of SPACE versus PACE, 100 sparse samples are created for each year. In each sample, 5 observations are randomly selected per curve. In the calculation of reconstruction error, the smoothed EVI based on dense observations that is used in the isotropy test serves as the underlying true process.}}
  			\label{fig: Harvard Forest isotropy and reconstruction}		
  		\end{figure}

		The smoothing method is described in \citet{liu2012functional} and more general introduction can be found in \citet{RamsaySilverman2005}. It has been shown that smoothing of individual curves and covariance surface are asymptotically equivalent. When applying the isotropy test proposed in Section \ref{subsection: isotropy test} to Harvard Forest EVI, we swap the estimation of fPC scores and anisotropy \matern parameters. In particular, we estimate fPCs and associated scores of hypothesized null curves through the process described above. Next, we estimate anisotropy \matern parameters based on the empirical correlation calculated from estimated fPC scores. Following the proposed procedure, we only resample the second fPC scores as they present potential anisotropy effect. For simplicity, we construct noiseless hypothesized curves so re-smoothing is not needed. The null distribution of the estimated anisotropy angle $\alpha$ is shown in Figure \ref{fig: Harvard Forest isotropy and reconstruction}. The result suggests the rejection of isotropy effect and confirms the diagonal pattern in fPC scores.
								
			Next, we apply SPACE framework to the gap-filling of Harvard Forest EVI data. To that end, we create sparse samples by randomly selecting 5 observations
from each location and each year. 100 sparse samples are created. To make the estimation and reconstruction of EVI across 625 pixels more computationally tractable, both SPACE and PACE are performed in each year respectively. The distribution of IP in year 2005 is summarized in Figure \ref{fig: Harvard Forest isotropy and reconstruction}. 100\% of the 100 samples show improved reconstruction performance using SPACE. Similar results are also achieved for other years. By incorporating spatial correlation estimated from EVI data, better gap-filling can be achieved. 
  		
  \section{Conclusion}\label{section: conclusion}    
  		Much of the literature in functional data analysis assumes no spatial correlation or ignoring spatial correlation if it is mild. We propose the spatial principal analysis based on conditional expectation (SPACE) to estimate spatial correlation of functional data, using non-parametric smoothers on curves and surfaces. We show that the leave-one-bin-out cross-validation based on binned data performs well in selecting bandwidth for local linear smoothers. Empirical spatial correlation is calculated as the ratio of eigenvalues of cross-covariance and covariance surfaces. 
  		
  		A parametric model, \matern correlation augmented with anisotropy parameters, is then fitted to empirical spatial correlations at a sequence of spatial separation vectors. With finite sample, estimates are better for separable covariance than non-separable covariance. The fitted anisotropy \matern parameters can be used to compute the spatial correlation at any given spatial separation vector and thus are used to reconstruct trajectories of sparsely observed curves. 
			
			This work compares with the work in \citet{YaoMullerWang2005} where curves are assumed to be independent. We show that by incorporating the spatial correlation, reconstruction performance is improved. It is observed that the higher the noise and true spatial correlations, the greater the improvements. We demonstrate the flexibility of SPACE model in modeling the separability and anisotropy effect of covariance structure. Moreover, two tests are proposed as well to explicitly answer if covariance is separable and/or isotropy. Resonable empirical sizes and powers are obtained in each test. 
			
			Then we apply the SPACE method to Harvard Forest EVI data. In particular, we confirm the diagonal pattern observed in the second fPC scores through a slightly modified version of the proposed isotropy test. Moreover, we demonstrate that by taking into account explicitly the spatial correlation, SPACE is able to produce more accurate gap-filled EVI trajectory on average.

			In the end, we also present a series of asymptotic results which demonstrate the consistency of model estimates. In particular, the same convergence rate for correlated case as that of i.i.d. case is derived assuming mild spatial correlation structure.

\section*{\Large Appendix}
\small

	\begin{appendices}	

	\section{More on Consistency}\label{appendix: consistency}
		This section discusses modifications and improvements of results in \citet{YaoMullerWang2005} in order to incorporate spatial correlation. Theorem 1* of Section \ref{section: consistency} extends the corresponding result of \citet{YaoMullerWang2005} and serves as the foundation of our work. In subsequent discussions, we change some notations in \citet{YaoMullerWang2005} to accommodate existing ones in our work. To avoid confusions, we make the following clarifications. We denote the number of curves by capital $N$ whereas Yao et.al. used $n$. For the number of observations on curve $i$, we use $n_i$ whereas Yao et.al. used $N_i$. Yao et.al. assumed $N_i$ follows the distribution of a random variable $N$ which is denoted by $N^*$ in our work. We denote all variables which depend on the BLUP of fPC scores by check sign instead of tilde used by Yao et.al. We use the same notations as those in \citet{YaoMullerWang2005} for other variables unless stated otherwise.  

		\citet{YaoMullerWang2005} states five theorems together with a series of lemmas and corollaries to describe the convergence of PACE estimators. Theorem 1 establishes the uniform convergence rates in probability of local linear estimates of mean $\mu(t)$ and covariance function $G(s,t)$. Corollary 1 states the convergence rate in probability of measurement error variance $\sigma^2$ and is derived by applying the results in Theorem 1. Theorem 2 establishes the convergence in probability of the eigen-structure of the covariance kernel, which includes the convergence of eigenvalues $\lambda_k$ of multiplicity 1, $L^2$ and uniform convergence of eigenfunctions $\phi_k(t)$'s associated with $\lambda_k$ of multiplicity 1. Based on the results in Theorem 1 and 2, derived in Theorem 3 are the convergence of estimated fPC score $\hat{\xi}_{ik}$ towards the conditional expectation of true fPC score (BLUP) $\check{\xi}_{ik}$, and the pointwise convergence of reconstructed curves $\widehat{X}_{K,i}(t)$ based on K-dimensional approximation towards its associated infinite dimensional target $\check{X}_i(t)$ constructed from $\{\check{\xi}_{ik}\}_{i=1,k=1}^{N,\infty}$. With Gaussian assumption, Theorem 4 establishes the pointwise asymptotic distribution of the standardized deviance between reconstructed curve $\widehat{X}_{K,i}(t)$ and true curve $X_i(t)$. In Theorem 5, simultaneous confidence region of $X_{K,i}(t)$ over $\cal T$ in terms of reconstructed curves $\widehat{X}_{K,i}(t)$ is established. Finally, the simultaneous region of K-dimensional vector ${\bm\xi}_{K,i}$ in terms of estimated fPC score vector $\widehat{\bm\xi}_{K,i}$  is derived in Corollary 2. As the cornerstone of all other derived results, Theorem 1 is derived based on the results in Lemma 1 and 2. These two lemmas establish the uniform convergence rates for the weighted average $\Psi_{pN}(t)$ of function $\psi_p(T_{ij}, Y_{ij})$ and kernel function $\kappa_1((t-T_{ij})/h_{\mu})$ for one-dimensional smoother, and the weighted average $\Theta_{pN}(t,s)$ of function $\theta_p(T_{ij},T_{il},Y_{ij},Y_{il})$ and kernel function $\kappa_2((t-T_{ij})/h_G,(s-T_{il})/h_G)$ for two-dimensional smoother. These weighted averages are the building blocks of the closed form expressions of the local linear smoothers. For example, the local constant estimator of $\mu_t$ can be expressed as 
		\begin{eqnarray}
			\dfrac{\sum_i\sum_j\kappa_1\left(\dfrac{t-T_{ij}}{h_{\mu}}\right)Y_{ij}}{\sum_i\sum_j\kappa_1\left(\dfrac{t-T_{ij}}{h_{\mu}}\right)} = \dfrac{\dfrac{1}{Nh_{\mu}}\sum_i\dfrac{1}{N^*}\sum_j\kappa_1\left(\dfrac{t-T_{ij}}{h_{\mu}}\right)Y_{ij}}{\dfrac{1}{Nh_{\mu}}\sum_i\dfrac{1}{N^*}\sum_j\kappa_1\left(\dfrac{t-T_{ij}}{h_{\mu}}\right)} = \dfrac{\Psi_{1N}}{\Psi_{2N}}.
		\end{eqnarray}
			
		The above results are only valid based on the assumption of a series of conditions. \citet{YaoMullerWang2005} describes three groups of conditions named after letters A,B and C. Some more general assumptions are also made. Specifically, both time point $T_{ij}$ and the observed value of curve $i$ made at this time point, $Y_{ij} = Y_i(T_{ij})$, are random variables. $T_{ij}$ is assumed to follow the same marginal distribution $T$ with density $f(t)$, data pair $(T_{ij}, Y_{ij})$ is assumed to follow the same joint distribution $(T,Y)$ with density $g(t,y)$, for any given curve $i$ and time index $j$. The random vector $(T_{ij},T_{il},Y_{ij},Y_{il})$ for $j\neq l$ and any $i$ is assumed to follow the same distribution $(T_1,T_2,Y_1,Y_2)$ with density $g(t_1,t_2,y_1,y_2)$. In addition, $T_{ij}$ is assumed to be independent across all $i$ and $j$, and both $(T_{ij}, Y_{ij})$ and $(T_{ij},T_{il},Y_{ij},Y_{il})$ are independent across curve index $i$. 
		
		Conditions series (A1) pertains to the number of observations per curve. It is assumed to be a discrete random variable following the distribution $N^*$. $N^*$ is assumed to have finite mean and takes value greater than 1 with positive probability, and is independent of all $T_{ij}$'s and $Y_{ij}$'s. Series (A2) specifies the convergence rates of smoothing bandwidths and their relative rates with respect to the number of curves $N$. This series of conditions ensures the convergence rate derived in Theorem 1. Series (A3) and (A4) assume favorable properties of smoothing kernel functions and finite fourth centered moments of $Y$, such that the auxiliary results in Lemma 1 and 2 hold. Condition (A5) assumes Gaussian distribution of fPC scores required by Theorem 4 and 5 and Corollary 2. Condition (A6) assumes that the data asymptotically follow a linear scheme. Condition (A7) assumes the pointwise limit of $\bm\phi_{K,s}\bm\Sigma(\check{\bm\xi}_{K,i}-\bm\xi_{K,i})\bm\phi_{K,t}\triangleq\bm\phi_{K,s}\bm\Omega_K\bm\phi_{K,t}$ as $K\rightarrow\infty$ exists such that Theorem 4 and 5 hold. Category B specifies properties of density and smoothing kernel functions for Lemma 1 and 2 to hold. Category C assumes favorable continuity and boundedness properties of functions $\psi_p$ and $\theta_p$ crucial for Lemma 1 and 2. 
		
		To reduce notational confusion between $v$ and $\nu$ in \citet{YaoMullerWang2005}, we denote the order of kernel functions $\kappa_1(t)$ and $\kappa_2(t,s)$ by $(a,b)$ instead of $(\nu,l)$ which is used in \citet{YaoMullerWang2005}. We review two lemmas and five theorems in more details below. The following review provides a skeleton of proofs in \citet{YaoMullerWang2005}, which serves as a scaffold to facilitate the introduction of our results. Please refer to \citet{YaoMullerWang2005} for more details regarding proofs for the i.i.d. case. 
		
		Proofs in \citet{YaoMullerWang2005} refer to Slutsky's theorem in many occasions. However, we found this reference is not precise because the statement of Slutsky's theorem as the widely known version does not involve convergence rate and big O notation as \citet{YaoMullerWang2005} have assumed and suggested. Since it is heavily used throughout the proofs, we state this auxiliary result in a separate lemma, Lemma 4, indexed after the existing lemmas in \citet{YaoMullerWang2005}.
			
			{\bf Lemma 4}
			{\it Assume $\sup_{t\in{\cal T}}|V_n(t) - V(t)| = O_p(c_n)$ where $c_n\rightarrow 0$ as $n\rightarrow\infty$, and $g(u)\in C^2$ on ${\cal R}^d$ for any $d>0$. Then 
			\begin{equation}\label{lemma4}
				\sup_{t\in{\cal T}}|g(V_n)-g(V)| = O_p(c_n).
			\end{equation}
			}			
			
			Proof of this lemma is deferred to Appendix \ref{appendix: proofs}. As a special case, a stronger rate can be derived for $g(u) = u_1\times u_2$. Specifically, assume $\sup_{t}|V_{1n}(t)-V_1(t)|=O_p(a_n)$,	$\sup_{t}|V_{2n}(t)-V_2(t)|=O_p(b_n)$ and $a_n/b_n\rightarrow 0$ as $n\rightarrow\infty$. Then it is easy to show that $\sup_{t}|V_{1n}(t)V_{2n}(t)-V_1(t)V_2(t)|=O_p(a_n)$.
	
	Recall that in \citet{YaoMullerWang2005} $(T_{ij},T_{ik},Y_{ij},Y_{ik})$ is assumed to follow the distribution of random vector $(T_1,T_2,Y_1,Y_2)$ with density function $g(t_1,t_2,y_1,y_2)$. Accordingly, for each spatial separation vector $\bm\Delta$ and any location pairs $(i,j)\in N(\bm\Delta)$, we assume $(T_{ik},T_{jl},Y_{ik},Y_{jl})$ follows the distribution of $(T_1,T_2,Y_1,Y_2)$ with density function $g_{\ScriptSmall{\bm\Delta}}(t_1,t_2,y_1,y_2)$. We assume the same regularity conditions on $g_{\ScriptSmall{\bm\Delta}}(t_1,t_2,y_1,y_2)$ as those in \citet{YaoMullerWang2005}. As mentioned earlier, Lemma 1 and 2 in \citet{YaoMullerWang2005} are the foundations and rely on the assumption of zero spatial correlation. Modifications due to the incorporation of spatial dependence start from these two lemmas.
	
	{\bf Lemma} ${\bf 1^*}$, as the counterpart of Lemma 1 in \citet{YaoMullerWang2005}, states the uniform convergence of weighted averages $\Psi_{pN}(t)$ for one-dimensional case. As pointed out above, the key step of achieving the stated rate is to show that $Var(\varphi_{pN}(u)) = O(1/N)$. With the presence of spatial correlation, we have,
			\begin{eqnarray}
				Var(\varphi_{pN}(u)) &=& Var\left(\dfrac{1}{N}\sum_{i=1}^N\dfrac{1}{EN^*}\sum_{j=1}^{n_i}e^{-iuT_{ij}}\psi_p(T_{ij},Y_{ij})\right)\nonumber\\
														 &=& \dfrac{1}{N^2}Var\left(\sum_{i=1}^N\dfrac{1}{EN^*}\sum_{j=1}^{n_i}e^{-iuT_{ij}}\psi_p(T_{ij},Y_{ij})\right)\nonumber\\
														 &=& \dfrac{1}{N}Var\left(\dfrac{1}{EN^*}\sum_{j=1}^{N^*}e^{-iuT_{j}}\psi_p(T_j,Y_j)\right)\nonumber\\ 
														 &+& \dfrac{1}{N^2}\sum_{1\leq i\neq j\leq N}Cov\left(\dfrac{1}{EN^*}\sum_{k=1}^{n_i}e^{-iuT_{ik}}\psi_p(T_{ik},Y_{ik})\right.,\nonumber\\
														 && \left.\hspace{90pt}\dfrac{1}{EN^*}\sum_{l=1}^{n_j}e^{-iuT_{jl}}\psi_p(T_{jl},Y_{jl}))\right)\nonumber\\
														 &\triangleq& \dfrac{1}{N}Var\left(\dfrac{1}{EN^*}\sum_{j=1}^{N^*}e^{-iuT_{j}}\psi_p(T_j,Y_j)\right) + \dfrac{1}{N^2}\sum_{1\leq i\neq j\leq N}Cov(Z_i,Z_j)\nonumber\\
														 &\leq& \dfrac{1}{N}E\left[\left(\dfrac{1}{EN^*}\sum_{j=1}^{N^*}e^{-iuT_j}\psi_p(T_j,Y_j)\right)^2\right] + \dfrac{1}{N^2}\sum_{1\leq i\neq j\leq N}Cov(Z_i,Z_j)\nonumber\\
														 &\leq& \dfrac{1}{N(EN^*)^2}E\left[\left(\sum_{j=1}^{N^*}e^{-i2uT_j}\right)\left(\sum_{j=1}^{N^*}\psi_p^2(T_j,Y_j)\right)\right]\nonumber\\
														 &+& \dfrac{1}{N^2}\sum_{1\leq i\neq j\leq N}Cov(Z_i,Z_j)\nonumber\\
														 &\leq& \dfrac{1}{N(EN^*)^2}E\left[N^*\left(\sum_{j=1}^{N^*}\psi_p^2(T_j,Y_j)\right)\right] + \dfrac{1}{N^2}\sum_{1\leq i\neq j\leq N}Cov(Z_i,Z_j)\nonumber\\
														 &=& \dfrac{1}{N(EN^*)^2}E\left[N^*\left(\sum_{j=1}^{N^*}E\psi_p^2(T_j,Y_j)\Big|N^*\right)\right] + \dfrac{1}{N^2}\sum_{1\leq i\neq j\leq N}Cov(Z_i,Z_j)\nonumber\\
														 &=& \dfrac{1}{N(EN^*)^2}E\left[N^{*2}E(\psi_p^2(T,Y))\right] + \dfrac{1}{N^2}\sum_{1\leq i\neq j\leq N}Cov(Z_i,Z_j)\nonumber\\
														 &=& \dfrac{EN^{*2}}{N(EN^*)^2}E(\psi_p^2(T,Y)) + \dfrac{1}{N^2}\sum_{1\leq i\neq j\leq N}Cov(Z_i,Z_j)\label{var of varphi}
			\end{eqnarray}
			The first term in the last line is $O(1/N)$ provided that $E(\psi_p^2(T,Y)) < \infty$ and $E(N^*)^2/(EN^*)^2$ is O(1) as $N\rightarrow\infty$. $\sum_{1\leq i\neq j\leq N}Cov(Z_i,Z_j)/N^2$ is the extra term if we introduce spatial correlations. Then requiring $\sum_{1\leq i\neq j\leq N}Cov(Z_i,Z_j)/N^2 = O(1/N)$ is a straight forward way for $Var(\varphi(u))$ to remain $O(1/N)$. Note there are $N(N-1)$ items in the summation of extra term, we would expect the correlation between $Z_i$ and $Z_j$ to decay fast enough as locations $i$ and $j$ are further apart. In a sum, to achieve a desired convergence rate, off-diagonal entries of spatial correlation matrices need to decrease at a proper pace when getting away from the diagonal. Within the current context, we introduce condition (D1) as follows:
	\begin{mydescription}{34.5pt}
				\item[(D1)] Let $Z_i = \sum_{k=1}^{n_i}e^{-iuT_{ik}}\psi_p(T_{ik},Y_{ik})/EN^*$ for $i = 1,\cdots,N$. Assume, 
					\begin{equation}
						\sum_{1\leq i\neq j\leq N}Cov(Z_i,Z_j)/N < \infty.
					\end{equation}
	\end{mydescription}
	
	{\bf Lemma} ${\bf 2^*}$ is the two-dimensional counterpart of {\bf Lemma} ${\bf 1^*}$. Note the local linear smoother of $G_{\ScriptSmall{\bm\Delta}}(s,t)$ is obtained by minimizing
			\begin{equation}\label{app 2d cross cov smoother}
				\sum_{(i,j) \in N(\ScriptSmallCap{\bm\Delta})}\sum_{k = 1}^{n_i}\sum_{l = 1}^{n_j}\kappa_2(\dfrac{t_{ik} - s}{h_G},\dfrac{t_{jl} - t}{h_G})(D_{ij}(t_{ik},t_{jl}) - \beta_0 - \beta_1(s - t_{ik}) - \beta_2(t - t_{jl}))^2
			\end{equation}
			Accordingly, we change the definition of weighted average $\Theta_{pN}(s,t)$ as
			\begin{eqnarray*}
				\Theta_{pN} &=& \Theta_{pN}(t,s)\\
										&=& \dfrac{1}{|N(\bm\Delta)|h_G^{|a|+2}}\sum_{(i,j) \in N(\ScriptSmallCap{\bm\Delta})}\sum_{k = 1}^{n_i}\sum_{l = 1}^{n_j}\theta_p(T_{ik},T_{jl},Y_{ik},Y_{jl})\kappa_2\left(\dfrac{t-T_{ik}}{h_G},\dfrac{s-T_{jl}}{h_G}\right)
			\end{eqnarray*}
			Note the outermost layer of summation is taken over all location pairs in $N(\bm\Delta)$. We thus change all $N$ that accompanied by $h_G$ with $|N(\bm\Delta)|$ in the lemmas, theorems, corollaries and conditions. In particular, affected conditions include (A2.2) and (C2.1b) in \citet{YaoMullerWang2005}. To achieve the same rate in Lemma 2 of \citet{YaoMullerWang2005}, we derive the following upper bound for the variance of $\gamma_{pN}(u,v)$.
			\begin{eqnarray*}
					&&  Var(\gamma_{pN}(u,v))\\
					&=& \dfrac{1}{N(\bm\Delta)}Var\left(\dfrac{1}{(EN^*)^2}\sum_{k=1}^{N^*}\sum_{l=1}^{N^*}e^{-(iuT_{k}+ivT_{l})}\theta_p(T_k,T_l,Y_k,Y_l)\right)\\
					&+& \dfrac{1}{|N(\bm\Delta)|^2}\sum_{\footnotesize{\substack{(i_1,j_1)\in N(\bm\Delta)\\(i_2,j_2)\in N(\bm\Delta)\\(i_1,j_1)\neq(i_2,j_2)}}}Cov\left(\dfrac{1}{(EN^*)^2}\sum_{k_1=1}^{n_{i_1}}\sum_{l_1=1}^{n_{j_1}}e^{-(iuT_{k_1}+ivT_{l_1})}\theta_p(T_{k_1},T_{l_1},Y_{k_1},Y_{l_1}),\right.\\
					&& \left.\hspace{132pt}\dfrac{1}{(EN^*)^2}\sum_{k_2=1}^{n_{i_2}}\sum_{l_2=1}^{n_{j_2}}e^{-(iuT_{k_2}+ivT_{l_2})}\theta_p(T_{k_2},T_{l_2},Y_{k_2},Y_{l_2})\right)\\
					&\triangleq& \dfrac{1}{N(\bm\Delta)}Var\left(\dfrac{1}{(EN^*)^2}\sum_{k=1}^{N^*}\sum_{l=1}^{N^*}e^{-(iuT_{k}+ivT_{l})}\theta_p(T_k,T_l,Y_k,Y_l)\right)\\
					&+& \dfrac{1}{|N(\bm\Delta)|^2}\sum_{\footnotesize{\substack{(i_1,j_1)\in N(\bm\Delta)\\(i_2,j_2)\in N(\bm\Delta)\\(i_1,j_1)\neq(i_2,j_2)}}}Cov\left(Q_{i_1,j_1},Q_{i_2,j_2}\right)
		 \end{eqnarray*}
		 Therefore, condition (D2) is given as
			\begin{mydescription}{34.5pt}
					\item[(D2)] Let $Q_{i,j} = \sum_{k=1}^{n_i}\sum_{l=1}^{n_j}e^{-(iuT_k+ivT_l)}\theta_p(T_k,T_l,Y_k,Y_l)$ for all $(i,j)\in N(\bm\Delta)$. Assume, 
					\begin{equation}
						\sum_{\footnotesize{\substack{(i_1,j_1)\in N(\bm\Delta)\\(i_2,j_2)\in N(\bm\Delta)\\(i_1,j_1)\neq(i_2,j_2)}}}Cov\left(Q_{i_1,j_1},Q_{i_2,j_2}\right)/|N(\bm\Delta)| < \infty.
					\end{equation}
			\end{mydescription}
			
		{\bf Theorem} ${\bf 1^*}$ adds conditions (D1) and (D2) to those listed in Theorem 1 in \citet{YaoMullerWang2005} to ensure same convergence rates. In particular, the uniform convergence rate of the cross covariance estimator is stated as 
		$$
			\sup_{t,s\in{\cal T}}|\widehat{G}_{\ScriptSmall{\bm\Delta}}(s,t) - G_{\ScriptSmall{\bm\Delta}}(s,t)| = O_p\left(\dfrac{1}{\sqrt{|N(\bm\Delta)|}h_G^2}\right).
		$$
		
		{\bf Theorem} ${\bf 2^*}$ establishes the consistency of empirical spatial correlation ${\hat\rho}^*_k(\bm\Delta)$. We restate the entire theorem as follows,
		
		{\it
		Under (A1.1)-(A4) and (B1.1)-(B2.2b) in \citet{YaoMullerWang2005} with $a=0,b=2$ in (B2.2a), $a=(0,0),b=2$ in (B2.2b), and (D1)-(D2) that we introduce,
			\begin{eqnarray}
				|{\hat\lambda}_k(\bm\Delta) - \lambda_k(\bm\Delta)| &=& O_p\left(\dfrac{1}{\sqrt{|N(\bm\Delta)|}h^2_G}\right),\label{lambda consistency}
			\end{eqnarray} 
			\begin{equation*}
				{\hat\rho}^*_k(\bm\Delta) = {\hat\lambda}_k(\bm\Delta) / {\hat\lambda}_k((0,0))\overset{p}{\rightarrow}\rho^*_k(\bm\Delta),
			\end{equation*}
			\begin{eqnarray}
				\|{\hat\phi}_k-\phi_k\|_H &=& O_p\left(\dfrac{1}{\sqrt{|N(\bm\Delta)|}h^2_G}\right),\label{phi l2 consistency}
			\end{eqnarray}
			and
			\begin{equation}
				\underset{t\in{\cal T}}{\sup}|{\hat\phi}_k(t)-\phi_k(t)| = O_p\left(\dfrac{1}{\sqrt{|N(\bm\Delta)|}h^2_G}\right)\label{phi uniform consistency}
			\end{equation}
		}
		Conditions series A, B and C are the same as those in Theorem 2 in \citet{YaoMullerWang2005}. The consistency of empirical spatial correlation follows by Slutsky's theorem.
			
		{\bf Theorem} ${\bf 3^*}$ also includes new conditions (D1) and (D2) in addition to the existing ones in Theorem 3 of \citet{YaoMullerWang2005}. The major difference lies in the proof. Specifically, we estimate all fPC scores $\widetilde{\bm\xi}$ in one conditional expectation as in equation \eqref{fPC score estimation}. For the $i$th curve and any integer $K\geq 1$, let ${\bf H}_K = \check{\widetilde{\bm\xi}} = \bm\Sigma(\widetilde{\bm\xi}) - \bm\Sigma(\check{\widetilde{\bm\xi}})$ and ${\bf H}_{K,ii}$ is the diagonal block corresponding to the $i$th curve. Let $w_K(s,t) = {\bm\phi}_{K,t}^{\ScriptSmall{T}}{\bf H}_{K,ii}{\bm\phi}_{K,t}$ for $t,s\in{\cal T}$ and ${\hat w}_K(s,t) = \widehat{\bm\phi}_{K,t}^{\ScriptSmall{T}}\widehat{\bf H}_{K,ii}\widehat{\bm\phi}_{K,t}$. Note the diagonal block of a positive definite matrix is also positive definite. Then $\{w_K(s,t)\}$ is a sequence of continuous positive definite functions. These modifications don't change the statement of condition (A7) in \citet{YaoMullerWang2005}.
		
	  Statements of {\bf Theorem} ${\bf 4^*}$ and {\bf Theorem} ${\bf 5^*}$ are the same as Theorem 4 and Theorem 5 in \citet{YaoMullerWang2005}, except for the addition of conditions (D1) and (D2).
		
		Finally, we restate {\bf Corollary} ${\bf 2^*}$ as follows,
			
			{\it
			Under the assumptions of Theorem $\it 5^*$,
			\begin{equation}
				\underset{N\rightarrow\infty}{\lim}P\left\{\underset{{{\textbf{\textrm l}}}\hspace{1pt}\in{\cal A}}{\sup}\dfrac{|{{\textbf{\textrm l}}}^{\ScriptSmall{\hspace{1pt}T}}(\widehat{\widetilde{\bm\xi}} - {\widetilde{\bm\xi}})|}{\sqrt{{{\textbf{\textrm l}}}^{\ScriptSmall{\hspace{1pt}T}}{\widehat{\bf H}}_K{{\textbf{\textrm l}}}}}\leq\sqrt{\chi^2_{d,1-\alpha}}\right\}\geq 1-\alpha,\label{fPC uniform consistency}
			\end{equation}		
		where $\chi^2_{d,1-\alpha}$ is the $(1-\alpha)$th percentile of the Chi-square distribution with $d$ degrees of freedom.
		}
		
		In this corollary, the simultaneous confidence region of all fPC scores are considered rather than that of fPC scores of each individual curve in \citet{YaoMullerWang2005}. Recall that $\widetilde{\bm\xi}$ is the vector of fPC scores from all curves stacked together and ${\bf H}_K$ is the covariance of $\check{\widetilde{\bm\xi}} - \widetilde{\bm\xi}$. The arguments to prove Corollary $\rm 2^*$ remain the same if we replace ${\bm\xi}_i$, $\widehat{\bm\xi}_i$, $\check{\bm\xi}_i$ and ${\widehat{\bm\Omega}}_K$ in \citet{YaoMullerWang2005} with $\widetilde{\bm\xi}$, $\widehat{\widetilde{\bm\xi}}$, $\check{\widetilde{\bm\xi}}$ and ${\widehat{\bf H}}_K$ respectively, observing that the validity of Corollary 2 in \citet{YaoMullerWang2005} depends only on the positive-definiteness of ${\widehat{\bm\Omega}}_K$ and joint Gaussian assumption on ${\bm\xi}_i$ which all apply to the counterparts of correlated case.   
	\section{Proof of Lemma 4}\label{appendix: proofs}
	\begin{proof}
		According to the assumption, for any $\epsilon > 0$, there exists $M_{\epsilon}$ such that 
		\begin{equation}\label{lemma 4 proof assumption}
			P\left\{\sup_{t\in{\cal T}}\dfrac{|V_n(t)-V(t)|}{c_n}>M_{\epsilon}\right\}<\epsilon
		\end{equation}
		where we assume $c_n$ is positive without loss of generality. For any given $N_0$ and $\epsilon$, we have
		\begin{eqnarray}
			&&P\left\{\dfrac{|g(V_n(t))-g(V(t))|}{c_n}>N_0, \forall t\in{\cal T}\right\}\\
			&=& P\left\{|g(V_n)-g(V)| > c_nN_0, \forall t\in{\cal T}\right\}\\
			&=& P\left\{|g(V_n)-g(V)| > c_nN_0\cap|V_n-V|\leq c_nM_{\epsilon/2}, \forall t\in{\cal T}\right\}\\
			&& + P\left\{|g(V_n)-g(V)| > c_nN_0\cap|V_n-V| > c_nM_{\epsilon/2}, \forall t\in{\cal T}\right\}\\
			&\triangleq& A + B
		\end{eqnarray}
		Under the assumption, we have
		$$
			B\leq P\left\{\dfrac{|V_n(t)-V(t)|}{c_n} > M_{\epsilon/2}\right\} < \epsilon / 2
		$$
		Choose $S$ as the smallest closed interval which contains
		$$
			\bigcup_{t\in{\cal T}}{[V(t)-\max(c_n)M_{\epsilon/2},V(t)+\max(c_n)M_{\epsilon/2}]}. 
		$$
		Note that $g'(u)$ is continuous. Let $C = \max_{u\in S}|g'(u)|$. Then $|C| < \infty$ and we have $|g(V_n(t))-g(V)|\leq|C||V_n(t) - V(t)|$. Thus,
		$$
			A\leq P\left\{|C||V_n(t)-V(t)|>c_nN_0\right\}
		$$
	  Choose $N_0 > |C|M_{\epsilon/2}$, then $A + B < \epsilon$. Observing that $|C|$ does not depend on $t$ and $n$. The result in \eqref{lemma4} follows.
		\end{proof}
		
		\section{Asymptotic Confidence Intervals and Regions}\label{appendix: intervals and regions}
		In this section, we derive asymptotic pointwise confidence bands for the underlying smooth curve $\widehat{X}_i(t)$ and asymptotic simultaneous confidence regions for fPC scores $\widetilde{\bm\xi}$. The covariance matrix of $\check{\widetilde{\bm\xi}}$ is $\bm\Sigma(\check{\widetilde{\bm\xi}}) = \bm\Sigma(\widetilde{\bm\xi},\widetilde{\bf Y})\bm\Sigma(\widetilde{\bf Y})^{-1}\bm\Sigma(\widetilde{\bf Y},\widetilde{\bm\xi})$. Note $\check{\widetilde{\bm\xi}} = E(\widetilde{\bm\xi}|\widetilde{\bf Y})$. Then we have $E(\check{\widetilde{\bm\xi}}\check{\widetilde{\bm\xi}}^{\ScriptSmall{T}}) = E(\check{\widetilde{\bm\xi}}\widetilde{\bm\xi}^{\ScriptSmall{T}})$. The risk of $\check{\widetilde{\bm\xi}}$ is measured as $\bm\Sigma(\check{\widetilde{\bm\xi}} - \widetilde{\bm\xi}) = \bm\Sigma(\widetilde{\bm\xi}) - \bm\Sigma(\check{\widetilde{\bm\xi}})\triangleq{\bf H}_K$. We write subscript $K$ next to existing symbols to indicate that the first $K$ eigenfunctions are used to approximate $X_i(t)$. In addition, we write $\bm\phi_{K,t}$ as the the first $K$ eigenfunctions evaluated at time $t$. With Gaussian assumptions, $(\check{\widetilde{\bm\xi}} - \widetilde{\bm\xi}) \sim {\cal N}({\bf 0}, {\bf H}_K)$. Let $\widehat{\bf H}_{K,ii}$ be the diagonal sub-block of $\widehat{\bf H}_K$ corresponding to curve $i$. Theorem ${\rm 4^*}$ which is the counterpart of Theorem 4 in \citet{YaoMullerWang2005} establishes that the distribution of $\widehat{X}_{K,i}(t) - X_i(t)$ may be asymptotically approximated by ${\cal N}(0, \widehat{\bm\phi}_{K,t}^{\ScriptSmall{T}}\widehat{\bf H}_{K,ii}\widehat{\bm\phi}_{K,t})$. Assuming ``weak'' spatial correlation described in conditions (D1) and (D2), we construct the asymptotic pointwise confidence intervals for $X_i(t)$,
		\begin{equation}
			\widehat{X}_{K,i}(t)\pm\Phi^{-1}\left(1-\dfrac{\alpha}{2}\right)\sqrt{\widehat{\bm\phi}_{K,t}^{\ScriptSmall{T}}\widehat{\bf H}_{K,ii}\widehat{\bm\phi}_{K,t})}
		\end{equation} 
		where $\Phi$ is the standard Gaussian cumulative distribution function. 
		
		Next, consider the construction of asymptotic simultaneous confidence bands. Let $X_{K,i}(t)=\sum_{k=1}^K\xi_{ik}\phi_k(t)$. Theorem ${\rm 5^*}$ which is the counterpart of Theorem 5 in \citet{YaoMullerWang2005} provides the asymptotic simultaneous band for $\{{\widehat X}_{K,i}(t)-X_{K,i}(t)\}$ for a given fixed $K$. The Karhunen-$\rm{Lo\grave{e}ve}$ theorem implies that $\sum_{t\in{\cal T}}E(X_{K,i}(t) - X_i(t))^2$ is small for fixed and sufficiently large $K$. Therefore, ignoring a remaining approximation error that may interpreted as a bias, we may construct $(1-\alpha)$ asymptotic simultaneous bands for $X_i(t)$ through
		\begin{equation}
			{\widehat X}_{K,i}(t)\pm\sqrt{\chi^2_{K,1-\alpha}\widehat{\bm\phi}_{K,t}^{\ScriptSmall{T}}\widehat{\bf H}_{K,ii}\widehat{\bm\phi}_{K,t})}
		\end{equation}
		where $\chi^2_{K,1-\alpha}$ is the 100$(1-\alpha)$th percentile of the chi-squared distribution with $K$ degrees of freedom. Because $\sqrt{\chi^2_{K,1-\alpha}}>\Phi^{-1}(1-\alpha/2)$ for all $K\geq 1$, the asymptotic simultaneous band is always wider than the corresponding asymptotic pointwise confidence intervals. 
		
		We then construct simultaneous intervals for all linear combinations of the fPC scores. Given fixed number $K$, let ${\cal A}\in{\cal R}^{NK}$ be any linear space with dimension $d\leq NK$. Then asymptotically, it follows from the uniform results in Corollary ${\rm 2^*}$, which is the counterpart of Corollary 2 in \citet{YaoMullerWang2005}, in Section \ref{section: consistency} that for all linear combination ${\textbf l}^{\ScriptSmall{\hspace{1pt}T}}\widetilde{\bm\xi}$ simultaneously, where ${\textbf l}\in{\cal A}$
		\begin{equation}
			{\textbf l}^{\ScriptSmall{\hspace{1pt}T}}\widetilde{\bm\xi}\in{\textbf l}^{\ScriptSmall{\hspace{1pt}T}}\widehat{\widetilde{\bm\xi}}\pm\sqrt{\chi^2_{d,1-\alpha}{\textbf l}^{\ScriptSmall{\hspace{1pt}T}}\widehat{\bf H}_K}{\textbf l}
		\end{equation}
		with approximate probability $1-\alpha$.
\end{appendices}

\bibliographystyle{Chicago}		
\bibliography{Dissertation}
\end{document}